\documentclass{amsart}
\usepackage[utf8]{inputenc}

\usepackage[english]{babel}

\usepackage{amssymb}
\usepackage{amsmath}
\usepackage{amsmath,amscd}
\usepackage{amsthm}
\usepackage{enumerate}
\usepackage{cite}
\usepackage{xcolor}
\usepackage[final]{listings}
\usepackage[section]{algorithm}

\makeindex

\newcommand{\C}{\mathbb{C}}
\newcommand{\Z}{\mathbb{Z}}
\newcommand{\Pj}{\mathbb{P}}

\newcommand{\N}{\mathbb{N}}

\newcommand{\G}{\mathbb{G}}

\newcommand{\vect}[1]{\mathbf{#1}}

\newcommand{\imm}{\operatorname{im}}

\newtheorem{thm0}{Theorem}[section]
\newtheorem{prop0}[thm0]{Proposition}

\newtheorem{coro0}[thm0]{Corollary}
\newtheorem{claim0}[thm0]{Claim}
\theoremstyle{definition}
\newtheorem{defn0}[thm0]{Definition}

\newtheorem{exa0}[thm0]{Example}

\newtheorem{rem0}[thm0]{Remark}


\author[E.~Angelini]{Elena Angelini}
\address{Dipartimento di Ingegneria dell'Informazione e Scienze Matematiche, Universit\`a di Siena, Italy}
\email{elena.angelini@unisi.it}

\author[L.~Chiantini]{Luca Chiantini}
\address{Dipartimento di Ingegneria dell'Informazione e Scienze Matematiche, Universit\`a di Siena, Italy}
\email{luca.chiantini@unisi.it}
\thanks{This work was supported by the National Group for Algebraic and Geometric Structures, and their Applications (GNSAGA – INdAM) and by the Italian PRIN 2015 - Geometry of Algebraic Varieties (B16J15002000005).}

\begin{document}

\title[Minimality and uniqueness for decompositions ]
{Minimality and uniqueness for decompositions of   specific ternary forms}
\date{}

\begin{abstract}
The paper deals with the computation of the rank and the identifiability of a specific ternary form. Often, one knows some short Waring decomposition of a given form, and the problem is to determine whether the decomposition is minimal and unique. We show how the analysis of the Hilbert-Burch matrix of 
the set of points representing the decomposition can solve this problem in the case of ternary forms. Moreover, when the decomposition
is not unique, we show how the procedure of liaison can provide alternative, maybe shorter, decompositions. We give an explicit 
algorithm that tests our criterion of minimality for the case of ternary forms of degree $9$. This is the first numerical case in which 
a new phaenomenon appears: the span of $18$ general powers of linear forms contains points of (subgeneric) rank $18$,  
but it also contains points whose rank is $17$, due to the existence of a second shorter decomposition which is completely different from the given one. 

\end{abstract}

\maketitle

\section{Introduction}

The paper deals with homogeneous polynomials (forms) $F$, which will be also seen as symmetric tensors, over the complex field $\C$.
Our analysis starts with a given {\it Waring }  expression of $F$ in terms of powers of linear forms 
$$ F= \lambda_1L_1^d+\dots +\lambda_rL_r^d$$
($d$ is the degree of $F$). Our target is to determine the minimality and the uniqueness of the expression, up to a scalar multiplication of the $L_i$'s. 
Indeed,  in our setting minimality is a consequence of uniqueness. So, when uniqueness holds, $r$ is the (Waring) rank of $F$.

The determination of the rank of a given form, as well as the the uniqueness of a Waring expression, 
is relevant in several aspects of tensor theory related to physics, signal processing, statistics, chemistry, artificial intelligence, etc.
For some examples of these relations, taken from the huge literature on the subject, let us mention the papers  \cite{AllmanMatiasRhodes09}, \cite{AnandkumarGeHsuKakadeTelgarsky14},\cite{AppellofDavidson81}, \cite{RaoLiZhang18}.

Our analysis, besides its theoretical interest, will produce the following concrete application. There are several situations in which one 
knows some Waring expression of a specific form $F$. For instance, this happens when $F$ splits in a sum 
of blocks whose rank is well known, or when one applies to $F$ some heuristic algorithm, e.g. based on local regression. In
these cases, however, the question about the minimality of the expression remains open. We will see in examples (see Section
\ref{exok}) that a local analysis of the Waring expression cannot guarantee its minimality, except when the length $r$
is rather small (as, e.g., in the Kruskal's criterion). We will produce a method which, in principle, can solve the minimality problem, 
and then guarantee that 
$r$ is the rank of $F$, for all Waring expressions of ternary forms. The consequent algorithm that one can construct
depends on the degree $d$. In the last section, we produce a concrete example of the algorithm, for ternary forms of degree
$9$.

We attack the problem with tools of  projective algebraic geometry. Thus, we start with the polynomial ring $R=
\C[x_0,\dots, x_n]$ and we indicate with $R_d$ the homogeneous piece of degree $d$. Any linear form $L$ determines
a point (which, by abuse of notation, we still indicate with $L$) in the projective space $\Pj^n$ 
over the linear space $R_1$ of forms of degree $1$. In this notation, $F$ corresponds to a point in the projective space over $R_d$.
We identify $\Pj(R_d)$, of (projective) dimension  $N=\binom{n+d} d -1$, with $\Pj(Sym^d(R_1))$. The map $\nu_d:\Pj^n\to \Pj^N$ which
sends a linear form $L$ to its $d$-th power is universally known as the $d$-th {\it Veronese} map.

A Waring expression $ F= \lambda_1L_1^d+\dots +\lambda_rL_r^d$ of $F$ determines a subset $A=\{L_1,\dots,L_r\}\subset\Pj^n$ such
that $F$ sits in the span of $\nu_d(A)$. In our notation, we will say that $A$ is a {\it decomposition} of length $r$ of $F$. The decomposition $A$ 
is {\it non-redundant} if $F$ does not lie in the span of $\nu_d(A')$, for any proper subset $A'\subset A$. The decomposition $A$ is 
{\it minimal} (resp. {\it unique}) if there are no other  subsets $B\subset \Pj^n$, of length $r'< r$ (resp. $r'\leq r$) such that $F$ belongs to the span of $\nu_d(B)$. 
So, when $A$ is minimal, its length is the (Waring) rank of $F$.

There are methods that, in some setting, can prove that a decomposition $A=\{L_1,\dots,L_r\}$ of $F$ is unique,
and also minimal. For instance,
one can take several flattenings of the tensors associated to $F,L_1^d,\dots, L_r^d $ and use them to prove the uniqueness, as in \cite{DomaLath13a}, \cite{DomaLath17}.
Alternatively, one can study the catalecticant map $\alpha$ defined by some partial derivatives of $F$, and the intersection of the image
of $\alpha$ with a suitable Veronese variety, as in \cite{LandOtt13},  \cite{MassaMellaStagliano18}. The most famous method to prove uniqueness goes back to
Kruskal, \cite{Kruskal77}. Kruskal used the {\it Kruskal's rank} of the matrix whose entries are the coefficients of the $L_i$'s (i.e. the coordinates
of the $L_i$'s, as points in $\Pj^n$) to produce an inequality which, if satisfied, guarantees that the decomposition is unique,
and minimal.
Kruskal's criterion has been refined in \cite{COttVan17b}, by taking into account the matrices defined by powers of the $L_i$'s.

All the aforementioned methods  have a common problem:  they will not provide an answer when the degree $d$ and
the length $r$ grow.

From the point of view of projective geometry, it is clear that the Waring rank is constant in a Zariski open subset of $\Pj^N$, i.e.
it is constant outside a subset defined by algebraic equations: a small subset, in any reasonable metric.
Denote with $r_g$ the {\it generic} value of the rank,  which holds in a dense subset of $\Pj^N$. Then, it is almost
straightforward that a decomposition cannot be unique when $r>r_g$ (too many parameters). When $r=r_g$ generic uniqueness can hold,
but it is a quite rare phenomenon, completely classified (see \cite{GaluppiMella}). On the contrary, for {\it subgeneric} values $r<r_g$, uniqueness 
holds for a minimal Waring expression of a {\it sufficiently general} form $F$, except for a short list of cases (see \cite{COttVan17a}). 
On the other hand, for a specific decomposition of a specific form $F$, whose rank is unknown, all the known methods can guarantee
its uniqueness or minimality only when the length $r$ is much smaller than $r_g$ (see Section 3 of \cite{AngeC}, for a concrete bound).

The analysis that we propose overruns this difficulty, and will guarantee the minimality, and also the uniqueness, of a given decomposition,
in principle for all $r<r_g$. We will give a concrete example for forms of degree $9$ in three variables.

The theoretical situation can be summarized as follows. Geometrically, the Waring expression $ F= \lambda_1L_1^d+\dots +\lambda_rL_r^d$ presents 
$F$ as a point of the (linear) span of the set $\nu_d(A)$. When $r$ is big, even if the $L_i$'s are general,
it turns out that the span of  $\nu_d(A)$ contains both points $F$ for which the decomposition $A$ is unique, and points 
 $F'= \mu_1L_1^d+\dots +\mu_rL_r^d$ for which uniqueness, and also minimality, do not hold  (though $A$ is still a non-redundant decomposition of $G$).
 Examples of this situation are described in Section \ref{sec:nonics18} below. It may happen indeed (Case 2 of Proposition \ref{possiblecases}) 
 that even if  $F'$ is non-redundantly spanned by $\nu_d(A)$, yet  there exists another decomposition $B$ of $F'$, completely different
 from $A$ and with less summands.
  Thus, methods based on the geometrical analysis of $A$ alone (as the Kruskal's or the catalecticant approaches) cannot 
  distinguish between forms $F,F'$ in the span  of $\nu_d(A)$, so they will fail to guarantee e.g. minimality, when the span 
  contains points with different behavior. The answer can be obtained only by analyzing $A$ {\it and} the coefficients $\lambda_i$'s of the Waring expression.
  
Our analysis follows the guidelines introduced in \cite{AngeC}. We attack the problem with a set of tools typical for the 
study of the geometry of finite sets in projective spaces: Hilbert-Burch matrices,  the Cayley-Bacharach
property, and liaison. An extensive illustration of the interplay between decompositions and the geometry of finite sets can be found in the book \cite{IK}. 
Papers \cite{BallC13}, \cite{BallBern12a}, \cite{BallBern13a} are based on the study of Hilbert functions of points. What is really new in the present paper, 
as well as in \cite{AngeC}, is the observation that given a finite set $A\subset \Pj^2$, the possible alternative (maybe even shorter) 
decompositions of any ternary form $F$ in the span of
$\nu_d(A)$ can be recovered geometrically, by playing with the Hilbert-Burch matrix of $A$ and liaison. Thus, we can characterize the (algebraic)
subset $\Theta$  of the span of $\nu_d(A)$ consisting of forms {\it for which $A$ is non-redundant}, but yet they have an alternative 
decomposition of length $r'< r$. Moreover,
we produce concrete algorithms which  guarantee that a form $F$ does not lie in the `bad' set $\Theta$.

In the paper, we illustrate in details the procedure in the case of ternary forms of degree $9$. For such forms, the generic rank
is $r_g=19$, see \cite{AlexHir95}. We consider the highest sub-generic value $r=18$. We show that for any generic choice of a subset $A$ of $18$ points
in $\Pj^2$ (the projective space of linear forms in three variables) the span $\Lambda$ of $\nu_9(A)$ contains: (i) points $F$
for which the decomposition $A$ is minimal and unique, (ii) points $F'$ for which $A$ in non-redundant, but there exists
an alternative decomposition $B\neq A$, of length $18$, (iii) points $F''$ for which $A$ is non-redundant, but there exists
an alternative decomposition of length $17$ (so their rank is smaller than $18$). Moreover, we prove that the (closure of) the
set of points satisfying (ii)  is a  hypersurface  of $\Lambda$, 
and it is a birational image of a Grassmannian of lines (see Theorem \ref{th9} below). We also determine properties of the closure $\Theta$
of the set of points satisfying (iii), which is a subvariety of $\Lambda$.

Finally, we produce an algorithm to test if  a given form $ F= \lambda_1L_1^9+\dots +\lambda_{18}L_{18}^9$ in the span $\Lambda$
of $\nu_9(\{L_1, \dots, L_{18}\})$ lies in the `bad' locus $\Theta$ defined above, i.e. its rank is smaller than $18$. 
The algorithm analyses both $A$ and the coefficients $\lambda_i$'s of the Waring expression. 
When the answer provided by the algorithm is negative, we can conclude that
$A$ is a minimal decomposition of $F$, so that $F$ has rank $18$. 
Similar algorithms for detecting the uniquess of the decomposition can be easily produced, though they need
more parameters. Finally, when a second decomposition of length $r'\leq r$ exists, we show how we can recover the new decomposition
from the known one $A$ (see the end of the examples in Section \ref{exok}).

The paper is structured as follows. In Sect. \ref{sec:notation} we introduce our setting and recall the main tools, coming from 
classical algebraic geometry, that characterize our analysis. In particular, we mention the Hilbert function and the Cayley-Bacharach 
property for finite sets in the projective space, the notion of Hilbert-Burch matrix and liason theory for sets of points in the projective plane. 
Sect. \ref{sec:nonics18} is the core of the paper: we develop our analysis for ternary forms of degree $ 9 $, dealing with ranks that are outside 
the range of applicability of the celebrated Kruskal's criterion. Sect. \ref{sec:intersection} disposes on the case in which two decompositions of a ternary nonic intersect. Finally, Sect. \ref{sec:alg} is devoted to the effective algorithm we developed 
according to the criterion of minimality (and uniqueness) obtained in Sect. \ref{sec:nonics18}. Several numerical examples are presented.
\smallskip

\paragraph{\textbf{Acknowledgements.} This work was supported by the National Group for Algebraic and Geometric Structures, 
and their Applications (GNSAGA – INdAM) and by the Italian PRIN 2015 - Geometry of Algebraic Varieties (B16J15002000005).}

\section{Preliminaries}\label{sec:notation}

\subsection{Notation}
For any finite set $A$ we denote with $\ell(A)$ the cardinality.\\
For $d,n \in \N$, let $ \C^{n+1} $ be the space of linear forms in $ x_{0}, \ldots, x_{n} $, so that $ S^{d} \C^{n+1}$ 
is the space of forms of degree $d$ in $n+1$ variables over $\C$. \\
Every form $ T \in S^{d} \C^{n+1} $ defines an element of $ \Pj(S^{d} \C^{n+1}) \cong \Pj^{N} $ $( N = \binom{n+d}{d} - 1) $, which, by abuse, we still denote by $T$. \\
We denote with $ \nu_{d}: \Pj^{n} \rightarrow \Pj^{N} $ the \emph{Veronese embedding} of $ \Pj^{n} $ of degree $ d $, which is given by 
$$ \nu_{d}([a_{0}x_{0}+ \ldots + a_{n}x_{n}]) = [(a_{0}x_{0}+ \ldots + a_{n}x_{n})^{d}]. $$ 
For any subset $Z\subset \Pj^N$ we denote with $\langle Z\rangle$ the linear span of $Z$. For instance, if $ A = \{P_{1}, \ldots, P_r\} \subset \Pj^{n} $ is a finite set,
$ \langle \nu_{d}(A) \rangle $ is the linear space in $\Pj^N$ spanned by the points $ \nu_{d}(P_{1}), \ldots, \nu_{d}(P_{r}) $. 
\smallskip

With the above notations we give the following definitions.

\begin{defn0}
Let $ A \subset \Pj^n $ be a finite set, and $ T \in S^{d} \C^{n+1} $ a form of degree $d$. 

We say that $ A $ \emph{computes} $ T $ if $ T\in \langle \nu_{d}(A) \rangle$.

We say that a set $ A $ which computes $T$ is \emph{non-redundant} if there are no proper subsets $ A' $ of $ A $ 
such that $A'$ computes $T$.

We say that a set $ A $ which computes $T$ is \emph{minimal} if there are no sets $B$, with $\ell(B)<\ell(A)$, 
such that $B$ computes $T$.

If $A$ computes $T$ and it is \emph{minimal}, the cardinality $\ell(A)$ is called the (Waring) \emph{rank} of $T$.
In this case we say that $A$ \emph{computes the rank} of $T$.
\end{defn0}

\begin{rem0}\label{rem:indep}
If $ A \subset \Pj^n $ is a finite set that computes the rank of $T$  then $A$ is non-redundant, and the points of $ \nu_{d}(A) $ are 
linearly independent, i.e.
$$ \dim(\langle\nu_{d}(A)\rangle) = \ell(A) -1.  $$
\end{rem0}

\begin{defn0}
A form $T$ is \emph{identifiable} if there exists a unique set $ A $ that computes the rank of $T$. 
\end{defn0}

\subsection{Kruskal's criterion for symmetric tensors}\label{sec:Kr}

One of the most used criteria for detecting the identifiability of a form $T$ is here described,
with some extensions, in geometric terms and adapted to the case of forms (= symmetric tensors).\\
We start with a definition which is the geometric analogue of the Kruskal's rank of a matrix.

\begin{defn0}
The \emph{d-th Kruskal's rank} of a finite set $ A \subset \mathbb{P}^{n} $ is 
$$ k_{d}(A) = \max \,\{k \, | \mbox{ for all } A' \subset A \mbox{ with } \ell(A') \leq k,\mbox{ then } \dim \langle\nu_{d}(A')\rangle = \ell(A') - 1 \}. $$
In other words, $k_d(A)$ is the maximal $k$ such that the image under $\nu_d$ of any subset of $A$ of cardinality
at most $k$ is linearly independent.
\end{defn0} 

\begin{rem0}\label{maxKrank}
For any $ d $, it holds that $ k_{d}(A) \leq \min\{N+1,\ell(A)\} $. 

Moreover,  if $k_d(A)=\min \{N+1,\ell(A)\}$ (i.e. $k_d(A)$ is maximal), then for all $A'\subset A$ the $d$-th Kruskal's rank
$k_d(A')$ is also maximal. 

If $ A $ is sufficiently general, then $ k_{d}(A) = \min\{N+1,\ell(A)\} $ i.e. the $d$-th Kruskal's rank is maximal (see e.g. Lemma 4.4 of \cite{COttVan17b}).
\end{rem0}

The Kruskal's rank is fundamental in the statement of the reshaped Kruskal's criterion.
 
\begin{thm0}[Reshaped Kruskal's Criterion, see \cite{COttVan17b}]\label{thm:kr}
Assume $ d \geq 3 $ and let $ A \subset \mathbb{P}^{n} $ be a non-redundant set computing $ T \in \Pj(S^{d} \C^{n+1}) $.  
Fix a partition  $ d = d_{1}+d_{2}+d_{3} $ with $ d_{1} \geq d_{2} \geq d_{3} \geq 1 $. If 
\begin{equation}\label{eq:Kr}
\ell(A) \leq \frac{k_{d_{1}}(A)+k_{d_{2}}(A)+k_{d_{3}}(A)-2}{2}
\end{equation}
then $ T $ has rank $ \ell(A) $ and it is identifiable.
\end{thm0}

In the case of ternary forms, the reshaped Kruskal's criterion has been recently extended in \cite{AngeC}, \cite{Ball19}, \cite{MourOneto}. 

\begin{thm0}\label{range} Let $ T \in \Pj(S^{d} \C^{3}) $ and let $ A \subset \mathbb{P}^{2} $ be a non-redundant set computing $ T $. 
Then $T$  is identifiable of rank $r$  if one of the following holds:

\begin{itemize}
\item $d=2m$ is even, $k_{m-1}(A)= \min\{\binom{m+1}2,r\}$,  $h_A(m)=r\leq \binom{m+2}2-2$;

\item $d=2m+1$ is odd, $k_{m}(A)=\min\{\binom{m+2}2,r\}$, $h_A(m+1)=r\leq \binom{m+2}2 +\lfloor {m \over 2} \rfloor$.

\end{itemize}
\end{thm0}

\begin{thm0}\label{ranger} Let $ T \in \Pj(S^{d} \C^{3}) $ and let $ A \subset \mathbb{P}^{2} $ be a non-redundant set computing $ T $. 
Then $A$ computes the rank of  $T$  if one of the following holds:

\begin{itemize}
\item $d=2m$ is even,  and $h_A(m)=r (\leq \binom{m+2}2)$;

\item $d=2m+1$ is odd, $k_{m}(A)=\min\{\binom{m+2}2,r\}$, $h_A(m+1)=r\leq \binom{m+2}2 + \lceil {m \over 2} \rceil $.

\end{itemize}
\end{thm0}

\subsection{The Hilbert function of finite sets and its difference}\label{sec:hilb}

\begin{defn0}
The \emph{evaluation map} of degree $ d $ on an ordered finite set of vectors $ Y = \{Y_{1}, \ldots, Y_{\ell}\} \subset \mathbb{C}^{n+1} $ is the linear map given by
$$ ev_{Y}(d): S^{d}\mathbb{C}^{n+1} \longrightarrow \mathbb{C}^{\ell} $$
$$ ev_{Y}(d)(F) = (F(Y_{1}), \ldots, F(Y_{\ell})). $$
\end{defn0}

\begin{defn0}\label{Hilbdef}
Let $ Y $ be a set of homogeneous coordinates for a finite set $ Z $ of $ \mathbb{P}^{n} $. The \emph{Hilbert function} of $ Z $ is the map
$$ h_{Z}: \mathbb{Z} \longrightarrow \mathbb{N} $$
such that $ h_{Z}(j) = 0$, for $j < 0$, $ h_{Z}(j) = rank (ev_{Y}(j))$, for $ j \geq 0. $
It is elementary indeed that the map depends only on $Z$, and not on the choice of a set of homogeneous coordinates for the points of $Z$.
\end{defn0} 

We will often use the \emph{first difference of the Hilbert function} $ Dh_{Z} $, given by
$$ Dh_{Z}(j) = h_{Z}(j)-h_{Z}(j-1), \,j \in \Z . $$

We collect some useful elementary properties of $ h_{Z} $ and $ Dh_{Z} $. They are well know and contained 
(sparsely) in the literature, see e.g. \cite{IK}, \cite{Migliore}. A summary can be found in \cite{C19}.

\begin{rem0}\label{rem:triv} 
\begin{enumerate}[(i)]
\item $ Dh_{Z}(j) \geq 0 $ for all $ j $; $ h_Z(j)=Dh_{Z}(j) = 0$ for $ j < 0$;  $ h_{Z}(0) = Dh_{Z}(0) = 1$.
\item $ h_{Z}(j) = \ell(Z) $ for all $ j \gg 0$, so that $ Dh_{Z}(j) = 0$ for $ j \gg 0 $ and $ \sum Dh_{Z}(j) = \ell(Z) $.
\item If $ Z' \subset Z $, then for all $ j \in \Z $  we have $ h_{Z'}(j) \leq h_{Z}(j) $ and $Dh_{Z'}(j) \leq Dh_{Z}(j).$
\end{enumerate}
\end{rem0}

\begin{prop0}\label{nonincr}
For any $  i >0 $ such that $ Dh_{Z}(i) \leq i $, it holds 
$$ Dh_{Z}(i) \geq Dh_{Z}(i+1), $$
i.e. the function $Dh_Z$ becomes non-increasing from $i$ on. 
Therefore, if $ Dh_{Z}(i) = 0 $, then $ Dh_{Z}(j) = 0 $ for any $ j \geq i $.
\end{prop0}

The following proposition is a straightforward application of the Grassmann formula in projective spaces.

\begin{prop0}\label{cap} 
Let $A,B\subset \Pj^n$ be disjoint finite sets such that both $\nu_d(A)$ and $\nu_d(B)$ are linearly independent. Set $Z=A\cup B$.\\
For any $d \in \N$, 
$$\dim(\langle v_d(A)\rangle\cap \langle v_d(B)\rangle) =\ell(Z)-h_Z(d)-1.$$
\end{prop0}

As pointed out e.g. in \cite{AngeC}, it follows that when $ A, B$ are two disjoint sets that compute a form $T$,
then $\langle v_d(A)\rangle\cap \langle v_d(B)\rangle$ is non empty, thus the union $Z=A\cup B$ satisfies $h_Z(d)<\ell(Z)$.
It follows:

\begin{prop0}\label{d+1} 
Let $ T \in S^{d}\C^{n+1} $ and let $ A, B \subset \Pj^{n} $ be disjoint non-redundant finite sets computing $ T $. 
Pose $ Z = A \cup B \subset \Pj^{n}$. Then $ Dh_{Z}(d+1) > 0 $.
\end{prop0}

The following Theorem has been generalized to sets of points in any projective space $\Pj^n$ (see \cite{BigaGerMig94}), but we
will need only the case $n=2$.

\begin{thm0}[Davis, \cite{Davis85}]\label{thm:Davis}
Let $ Z \subset \Pj^{2} $ be a finite set. Assume that for some $j>0$ one has
$ 0<Dh_{Z}(j) = Dh_{Z}(j+1)\leq j  $.\\
Then $Z$ splits in a union $ Z = Z_{1} \cup Z_{2} $ where $ Z_{1} $ lies on a curve of degree $ e=Dh_{Z}(j)$ and 
 $$ Dh_{Z_1}(i) = Dh_{Z}(i)  \mbox{ for } i\geq j.$$
\end{thm0}

The consequences of Davis' Theorem in our analysis are resumed in the following Proposition (see \cite{AngeC} Proposition 2.20).

\begin{prop0}\label{Dav}
Let $ T \in S^{d}\C^{3} $ and let $ A,B \subset \Pj^{2} $ be disjoint, non-redundant finite sets computing $ T $. Pose $ Z = A \cup B$.
Then there are no  indices $j\leq d$ such that  $ 0<Dh_{Z}(j) = Dh_{Z}(j+1) < j  $.
\end{prop0}
\begin{proof} (sketch) Assume on the contrary that $j$ exists. Then, by Davis' Theorem and the Grassmann formula, there exists a proper subset $Z_1\subset Z$
such that, if $A_1=A\cap Z_1$ and $B_1=B\cap Z_1$, then 
$$ \langle v_d(A)\rangle\cap \langle v_d(B)\rangle = \langle v_d(A_1)\rangle\cap \langle v_d(B_1) \rangle.$$
Thus, at least one between $A,B$ cannot be non-redundant.
\end{proof}

\subsection{The Cayley-Bacarach property}\label{sec:CB}

\begin{defn0}\label{def:CB}
A finite set $Z\subset \Pj^n$ satisfies the \emph{Cayley-Bacharach property in degree $i$}, abbreviated as $\mathit{CB}(i)$, if for all $P \in Z$, 
it holds that every form of degree $i$ vanishing at $ Z\setminus\{ P\}$ also vanishes at $P$.
\end{defn0}

The main consequence of the Cayley-Bacharach property on the Hilbert function of a set $Z$ is contained in the following two results,
proved in \cite{AngeCVan18}, by means of a deep result in \cite{GerKreuzerRobbiano93}.

\begin{thm0}\label{GKRext}
If a finite set $ Z \subset \Pj^{n} $ satisfies $\mathit{CB}(i)$, then for any $ j $ such that $ 0 \leq j \leq i+1 $ we have
$$ Dh_{Z}(0)+Dh_{Z}(1)+\cdots + Dh_{Z}(j) \leq Dh_{Z}(i+1-j)+\cdots +Dh_{Z}(i+1).$$
\end{thm0}

\begin{coro0}\label{CBconseq}
Let $ T \in S^{d}\C^{3} $ and let $ A \subset \Pj^{2} $ be a non-redundant finite set computing $ T $. 
Let $ B \subset \Pj^{2} $ be another non-redundant finite set computing $ T $
and assume $ A \cap B = \emptyset $. Then $Z=A\cup B$ satisfies the Cayley-Bacharach property $\mathit{CB}(d)$.
\end{coro0}

\subsection{The Hilbert-Burch matrix of a set of points in $\Pj^2$}\label{HBsec}

Let $A\subset \Pj^2$ be a finite set and call $I_A$ the homogeneous ideal of $A$, in the polynomial ring $R=\C[x_0,x_1,x_2]$. Then there exists 
an exact sequence of graded $R$-modules (called a \emph{minimal resolution} for $I_A$), as follows:
\begin{equation}
0 \longrightarrow F_1\xrightarrow M  F_0 \longrightarrow I_A \longrightarrow 0
\end{equation}
where $F_1,F_0$ are free graded $R$-modules, i.e. we have
$$   F_0 = \oplus_{i=1}^s R(-d_i)  \qquad  F_1 = \oplus_{j=1}^{s-1} R(-e_j)  $$
($R(i)$ is the graded module $R$ with the degrees shifted by $i$. The shift is needed in order to get that the maps in the
sequence are graded homomorphisms).\\
The map $F_0 \longrightarrow I_A $ sends the standard basis of $F_0$ to a set of minimal generators for $I_A$, and
the numbers $d_i$'s correspond to the degrees of the chosen minimal generators.\\
The map $F_1\longrightarrow  F_0$ is given by a matrix of forms $M$, which describes the \emph{first syzygies}
of $A$, i.e. the relations between the minimal generators.

\begin{defn0} \label{HBdef} The matrix of forms $M$ which determines the map $F_1\longrightarrow  F_0$ is called a \emph{Hilbert-Burch matrix}
for $A$. \\
The degree matrix of $M$ is called the  \emph{degree Hilbert-Burch matrix} of $A$.
\end{defn0}

\begin{thm0} \label{HBthm} (Hilbert-Burch Theorem, see \cite{CiGerOre88}) The Hilbert-Burch matrix $M$ depends on a choice of a set of 
minimal generators for the ideal $I_A$, but  the degree Hilbert-Burch matrix depends only on $A$.\\
The minimal generators which determine the surjection $F_0 \longrightarrow I_A$ are precisely the $(s-1)\times (s-1)$
minors of the matrix $M$ (taken with the corresponding sign). 
\end{thm0}

\begin{exa0}
A set of points $Z\subset\Pj^2$ is \emph{complete intersection} if there are two plane curves $F,G$ such that $Z=F\cap G$. In this case, the homogeneous ideal 
$I_{Z}$ is generated by the forms $F,G$ which define the two curves. If $d_1,d_2$ denote the degrees of $F,G$ respectively, then a minimal resolution
of $I_{Z}$ looks like
$$0 \longrightarrow R(-d_1-d_2)\xrightarrow M  R(-d_1)\oplus R(-d_2) \longrightarrow I_Z \longrightarrow 0$$
and the Hilbert-Burch matrix $ M$ is given by $M=\begin{pmatrix} G \\ -F \end{pmatrix}$.
\end{exa0}

\subsection{Linked sets and mapping cone}\label{sec:link}

Let $A,B\subset\Pj^2$ be two sets of points such that $Z=A\cup B$ is complete intersection of two curves $F$ and $G$, of degrees respectively $d_1$ and $d_2$. 
In this case, we say that $A$ and $B$ are \emph{linked} by (or also that $B$ is the \emph{residue} of $A$ in) a complete intersection of type $(d_1,d_2)$.

When the two sets are linked,  a Hilbert-Burch matrix of $B$ can be found from a Hilbert-Burch matrix of $A$ via the \emph{mapping cone} procedure (for more details, we refer to \cite{Ferrand75}, \cite{PeskineSzpiro74}, and Proposition 5.2.10 of \cite{Migliore}).

The homogeneous ideal $I_B$ is given by 
$$I_B=I_Z:I_A=\{ f\in R: fI_{A}\subset I_Z\}.$$
The inclusion $I_Z\subset I_A$ determines on the resolutions of $I_A$ and $I_Z$ a commutative diagram:
\begin{equation}\label{mc}
\begin{matrix}
0 & \longrightarrow & R(-d_1-d_2) & \xrightarrow  {M'} &  R(-d_1)\oplus R(-d_2) & \longrightarrow & I_Z &\longrightarrow & 0  \\
   &                            & \downarrow \phi' &                        & \downarrow    \phi         &                            &  \downarrow   &        &      \\
0 &\longrightarrow  & F_1               & \xrightarrow M     & F_0                                 & \longrightarrow & I_A &\longrightarrow & 0
\end{matrix}
\end{equation} 
where the rightmost vertical map is the inclusion.

Then $I_B$ is the image of the dual map $  R(d_1)\oplus R(d_2) \oplus F_1^\vee \xrightarrow{(\phi', M')^\vee} R(d_1+d_2)$, twisted by $-d_1-d_2$. A resolution of $I_B$ is given by:
$$0 \to F_0^\vee(-d_1-d_2) \xrightarrow {(\phi,M)^\vee} R(-d_2)\oplus R(-d_1) \oplus F_1^\vee(-d_1-d_2) \xrightarrow {(\phi', M')^\vee} I_B \to 0.$$

Notice that the resolution of $I_B$ obtained by the mapping cone procedure needs not to be minimal, in the sense that some summands of the map
$F_0^\vee(-d_1-d_2) \to R(-d_2)\oplus R(-d_1) \oplus F_1^\vee(-d_1-d_2) $ could be factored out.

\begin{exa0} Let $A$ be a general set of three points in $\Pj^2$. The ideal $I_A$ is generated by three quadrics, and a minimal resolution is given by
 $$0 \longrightarrow R(-3)^2\xrightarrow M  R(-2)^3 \longrightarrow I_A \longrightarrow 0,$$
 where the Hilbert-Burch matrix $M$ is a $3\times 2$ matrix of linear forms $\ell_{ij}$.
 If we take a general quadric $F$ and a general cubic $G$ containing $A$, we get a linkage between $A$ and another set $B$ of three points in the plane.
 Diagram \eqref{mc} looks like
 $$ \begin{matrix}
  0 & \longrightarrow & R(-5) & \longrightarrow  &  R(-2)\oplus R(-3) & \longrightarrow & I_Z &\longrightarrow & 0  \\
   &                            & \downarrow \phi' &                        & \downarrow    \phi         &                            &  \downarrow   &        &      \\
0 &\longrightarrow  & R(-3)^2               & \xrightarrow M     & R(-2)^3                      & \longrightarrow & I_A &\longrightarrow & 0
\end{matrix} $$
and the mapping cone gives the resolution
$$0 \to R(-3)^3 \xrightarrow {(\phi, M)^\vee} R(-3) \oplus R(-2) \oplus R(-2)^2  \longrightarrow  I_B \to 0.$$
The matrix of $(\phi,M)^\vee$ is obtained as follows. Since the minimal generators of $I_A$ are the minors of $M$, $F,G$ have a
representation as determinants of matrices as follows
$$ F=\det\begin{pmatrix} \ell_{11} & \ell_{12} & c_{13} \\  \ell_{21} & \ell_{22} & c_{23} \\ \ell_{31} & \ell_{32} & c_{33}\end{pmatrix} \quad
G=\det\begin{pmatrix} \ell_{11} & \ell_{12} & h_{13} \\  \ell_{21} &\ell_{22} & h_{23} \\  \ell_{31} & \ell_{32} & h_{33}\end{pmatrix}
$$
for some choice of constants $c_{i3}$ and linear forms $h_{i3}$.
Thus the matrix of $(\phi,M)^\vee$ is given by
$$ \begin{pmatrix} c_{13} & c_{23} & c_{33} \\ h_{13} & h_{23} & h_{33} \\ \ell_{11} & \ell_{21} & \ell_{31} \\  \ell_{12} & \ell_{22} & \ell_{32} 
\end{pmatrix}. $$
This resolution is not minimal, because we can choose the resolution of $I_A$ so that $F$ is the first generator. In this case
the row $(c_{13} , c_{23} , c_{33})$ is equal to $(1,0,0)$ and can be factored out. So, we can drop
the cubic generator of $I_B$. The resolution becomes
$$0 \to R(-3)^2 \xrightarrow {N} R(-2)\oplus  R(-2)^2 \longrightarrow  I_B \to 0$$
and the Hilbert-Burch matrix $N$ is:
$$ N= \begin{pmatrix}  h_{23} & h_{33} \\ \ell_{21} &  \ell_{31} \\   \ell_{22} &  \ell_{32} 
\end{pmatrix}. $$
\end{exa0}

\subsection{Grassmannian Varieties}

Let $ r,s \in \N $ such that $ r \leq s $ and let $ \G(r,s) $ be the Grassmannian variety of $r$-dimensional linear  spaces in $ \C^{s} $. $ \G(r,s) $ 
is an algebraic subvariety of a projective space by means of the Pl\"ucker embedding
\begin{equation}\label{eq:Pl}
i : \G(r,s) \hookrightarrow \Pj(\Lambda^{r}\C^{s}) 
\end{equation} 
defined by 
$$ i(\langle v_{1}, \ldots, v_{r} \rangle) = [v_{1} \wedge \ldots \wedge v_{r}]  $$
where $ v_{1}, \ldots, v_{r} $ are $ r $ linearly independent vectors of $ \C^{s} $. The homogeneous coordinates of 
$ \G(r,s) $ on $ \Pj^{\binom{s}{r}-1} \cong \Pj(\Lambda^{r}\C^{s}) $ are called Pl\"ucker coordinates. 
We recall that 
$$ \dim \G(r,s) = r(s-r). $$
Let $ W \in \G(r,s) $. We can associate to $ W $ the $ r \times s $ matrix with complex entries $ M_{W} $ whose rows contain the coordinates of a basis 
$ \{v_{1}, \ldots, v_{r}\} $ of $ W $. We notice that this representation is not unique: if we multiply for an element $U\in GL(r,\C) $, the new matrix represents the same point of $\G(r,s) $.
The Pl\"ucker coordinates of $ \G(r,s) $ are the minors of size $ r $ of $ M_{W} $, which are simply multiplied by a number when we substitute  $ M_{W} $ with $U M_{W} $.  \\

The embedding $i : \G(r,s) \hookrightarrow \Pj(\Lambda^{r}\C^{s}) $ is defined by a divisor $ h_{r,s}$ of $\G(r,s)$.
According to Lemma 11.1 of \cite{Dolgachev},  the divisor $ h_{r,s} \subset \G(r,s) $ associated to $ i $ is
$$ h_{r,s} = \{\Lambda \subset \C^{s} | \Lambda \cap \tilde{L} \not= \emptyset \}  $$
where $ \tilde{L} \cong \C^{s-r-1} $ is a fixed subspace of $\C^s$.

\section{Forms of degree $9$ in three variables}\label{sec:nonics18}

In this section, let us  assume that $ n = 2, d = 9 $ and let $ T \in S^{9} \C^{3} $. Thus $T$ is a form of degree $9$ in three variables, which
is associated to a curve of degree $9$ in $\Pj^2$.\\

Let $ A = \{P_{1},\ldots,P_{r}\} \subset \Pj^{2} $ be a finite set that computes $ T $.\\
Our target is to determine conditions on $T$ such that $A$ is the unique set
that computes the rank of $T$.\\ 

We will assume in this section that $A$ is sufficiently general, in a very precise sense.
 We assume indeed that:
\begin{itemize}
\item[(i)] $ A $ is non-redundant;
\item[(ii)] $ k_{4}(A) = \min\{15,r\}  $;
\item[(iii)] $ h_{A}(5) = r $.
\end{itemize}
It is a standard fact indeed that when $ r\leq 21$, all $A$ in a Zariski open subset of $(\Pj^2)^r$ satisfy the previous conditions.

\begin{rem0} For a general $ T \in S^{9} \C^{3} $, according to the Alexander-Hirschowitz Theorem \cite{AlexHir95}, the rank is $19$ and identifiability does not hold.\\
On the other hand, if $ r \leq 17 $ and $A$ satisfies the previous conditions, then, by Theorem \ref{range}, $ T $ is identifiable of rank $ r $. 
\end{rem0}

Therefore, we will assume in this section that $ r = 18 $ (so that $h_A(5)=18$ and $k_4(A)=15$).

Since $ A $ verifies properties $ (ii) $ and $ (iii) $, then the Hilbert function of $ A $ and its first difference verify, respectively,
\begin{equation}\label{eq:hDhA}
\begin{tabular}{c|cccccccc}
$j$ & $0$ & $1$ & $2$ &   $3$ &   $4$ &   $5$ &    $6$ &  $\dots$ \\  \hline
$h_{A}(j)$ & $1$ & $3$ &   $6$ &   $10$ &   $15$ &   $18$ & $ 18 $ & $\dots$ \cr
$Dh_{A}(j)$ & $1$ & $2$ &   $3$ &   $4$ &   $5$ & $3$ & $ 0 $ & $\dots$ \cr
\end{tabular}.
\end{equation}

\begin{rem0}\label{generic}The previous values of the Hilbert function imply that $A$ is not contained in quartic curves (it is clear, since we are assuming $k_4(A)=15$),
moreover $A$ is contained in $3$ independent quintics  $Q_{1}, Q_{2}, Q_{3}\in S^{5} \C^{3} $.
The quintics $Q_i$'s are elements of minimal degree in the homogeneous ideal $I_A$ of $A$. Thus, they are among the minimal generators.
If we multiply the $Q_i$'s by $3$ independent linear forms, which generate $S^{1} \C^{3}$, we obtain a set of nine forms of degree $6$ 
contained in the ideal $I_A$. These sextics span a subspace $\Lambda_6$ of dimension exactly nine in the homogeneous piece of
 degree $6$ of $I_A$. Indeed, the condition $ h_{A}(5) = r =18 $ implies $h_{A}(6)=18$, so that $(I_A)_6$, which is the kernel of the
 evaluation map in degree $6$, has affine dimension $28-18=10$, hence projective dimension $9$. 
 \end{rem0}
 
In section \ref{sec:alg}) we produce an algorithm which, starting with the coordinates of the points of $A$, texts
whether or not $A$ satisfies the generality conditions (i) - (iii).\\
Conditions (i) - (iii) determine the shape of the degree Hilbert-Burch matrix of $A$.

\begin{prop0}\label{HBA} Assume $A$ satisfies conditions (i) - (iii). Then there exist three forms $Q_{1}, Q_{2}, Q_{3}\in S^{5} \C^{3} $ of degree $5$ 
and one form $ S \in S^{6} \C^{3} $ of degree $6$ such that the ideal $I_A$ is minimally generated by $Q_1,Q_2,Q_3,S$. The resolution
of $I_A$, determined by this choice of generators, has the form
$$ 0 \longrightarrow R(-7)^{\oplus 3} \xrightarrow M  R(-5)^{3} \oplus R(-6) \longrightarrow I_A \longrightarrow 0.$$
\end{prop0}
\begin{proof} Since the space $(I_A)_5$ of quintics through $A$ is three-dimensional, a basis $Q_{1}, Q_{2}, Q_{3}$ for $(I_A)_5$ determines three
minimal generators for $I_A$. Since $ h_{A}(5) =18$, by the Castelnuovo-Mumford regularity Theorem (see \cite{Castelnuovo93}) $I_A$ is generated in degree $6$.
By Remark \ref{generic}, the three quintics $Q_{1}, Q_{2}, Q_{3}$ determine a subspace $\Lambda_6$ 
 of dimension $9$ in $(I_A)_6$. $(I_A)_6$ has dimension $10$, since $h_A(6)=18$.
It follows that by taking  $S\in (I_A)_6 \setminus \Lambda_6$, we get a minimal set of $4$ generators for $I_A$, which thus determines a surjection
$\alpha: R(-5)^{3} \oplus R(-6) \longrightarrow I_A$. The graded piece $(I_A)_7$ has dimension $\dim(R_7)-h_A(7)=36-18=18$. By multiplying each
$Q_i$ by a basis for $S^{2} \C^{3}$ and $S$ by a basis for $S^{1} \C^{3}$, we get $21$ elements in $(I_A)_7$. Thus we have at least $3$
independent relations of degree $7$ among $Q_{1}, Q_{2}, Q_{3},S$. Since $3=4-1$, by the Hilbert-Burch Theorem \ref{HBthm}, the independent relations
of degree $7$ are exactly $3$ and the kernel of the map $\alpha$ is $R^3(-7)$.
\end{proof}

By the Hilbert-Burch Theorem and by Proposition \ref{HBA},  there exist $ c_{uv}\in S^{2} \C^{3} $ and $ \ell_{j} \in S^1\C^{3} $, for $ u, v,j \in \{1,2,3\} $, 
such that $A$ has a Hilbert-Burch matrix $M$:
\begin{equation}\label{M}
M = \begin{pmatrix} c_{11} & c_{12} & c_{13}\\ c_{21} & c_{22} & c_{23}\\ c_{31} & c_{32} & c_{33}\\ \ell_{1} & \ell_{2} & \ell_{3} \end{pmatrix}.
\end{equation}
$ Q_{1}, Q_{2}, Q_{3}, S $ coincide, respectively, with $ (-1)^{i} $ times the minor obtained by leaving out the $ i $-th row of $M$, $ i \in \{1,2,3,4\} $.
In other words, we have:
\begin{multline}\label{geners}
Q_{1}=-\begin{vmatrix} c_{21} & c_{22} & c_{23}\\ c_{31} & c_{32} & c_{33}\\ \ell_{1} & \ell_{2} & \ell_{3}\end{vmatrix}, \quad
Q_{2}= \begin{vmatrix} c_{11} & c_{12} & c_{13}\\ c_{31} & c_{32} & c_{33}\\ \ell_{1} & \ell_{2} & \ell_{3}\end{vmatrix}, \\
Q_{3}=-\begin{vmatrix} c_{11} & c_{12} & c_{13}\\ c_{21} & c_{22} & c_{23}\\ \ell_{1} & \ell_{2} & \ell_{3}\end{vmatrix}, \quad
S = \begin{vmatrix} c_{11} & c_{12} & c_{13}\\ c_{21} & c_{22} & c_{23}\\ c_{31} & c_{32} & c_{33}\end{vmatrix}.
\end{multline}

\indent Now, assume that $ B = \{P'_{1},\ldots,P'_{\ell(B)}\} \subset \Pj^{2} $ is another finite set computing $ T $ such that
\begin{itemize}
\item[(i)] $\ell(B) \leq 18$;
\item[(ii)] $ B $ is non-redundant
\end{itemize}
and set $ Z = A \cup B \subset \Pj^{2} $.\\

{\it We assume, for the rest of the section, that the intersection $A\cap B$ is empty.}
\medskip

We will analyze the case in which the intersection is non-empty in the next section. 
Observe that, by Corollary \ref{CBconseq}, the last assumption:\\
 {\it implies that $Z$ satisfies the Cayley-Bacharach property $CB(9)$.}

\begin{rem0}
Since $ Z $ satisfies $ CB(9) $, then, by Theorem \ref{GKRext}, Remark \ref{rem:triv} (iii) and \eqref{eq:hDhA}, we get that
$$ Dh_{Z}(6) + Dh_{Z}(7) + Dh_{Z}(8) + Dh_{Z}(9) + Dh_{Z}(10) \geq $$
$$ \quad\quad\quad \geq  Dh_{A}(0) + Dh_{A}(1) + Dh_{A}(2) + Dh_{A}(3) + Dh_{A}(4) = 15. $$
Moreover, being $ Dh_{A}(5) = 3 $, then by  Remark \ref{rem:triv} (iii) we obtain that
\begin{equation}\label{eq:lbound}
36 \geq \ell(Z) \geq 15+ Dh_{Z}(5) + Dh_{Z}(6) + Dh_{Z}(7) + Dh_{Z}(8) + Dh_{Z}(9) + Dh_{Z}(10)  
\end{equation}
$$ \geq 18+ Dh_{Z}(6) + Dh_{Z}(7) + Dh_{Z}(8) + Dh_{Z}(9) + Dh_{Z}(10). \quad\quad\quad\, $$
Therefore 
\begin{equation}\label{eq:upbound}
Dh_{Z}(6) + Dh_{Z}(7) + Dh_{Z}(8) + Dh_{Z}(9) + Dh_{Z}(10) \leq 18. 
\end{equation}
We immediately get the following chain of inequalities
\begin{equation}\label{eq:chain}
15 \leq Dh_{Z}(6) + Dh_{Z}(7) + Dh_{Z}(8) + Dh_{Z}(9) + Dh_{Z}(10) \leq 18. 
\end{equation}
\end{rem0}

\begin{prop0}\label{possiblecases}
Assume that $ A \cap B = \emptyset $. Then one the following cases occurs for the first difference of the Hilbert function of $ Z $: \\
\textbf{Case 1}:
\begin{center}\begin{tabular}{c|cccccccccccccc}
$j$ & $0$ & $1$ & $2$ &   $3$ &   $4$ &   $5$ &  $6$ &  $7$ &  $8$ &  $9$ & $ 10 $ & $ 11 $ & $ \dots $ \\  \hline
$Dh_Z(j)$ & $1$ & $2$ &   $3$ &   $4$ &   $5$ &  $6$ &  $5$ &  $4$ & $3$ & $2$ & $ 1 $ & $ 0 $ & $ \cdots $ \cr
\end{tabular} \end{center}
and $ \ell(B) = 18 $, $ \ell(Z) = 36 $; \\
\textbf{Case 2}: 
\begin{center}\begin{tabular}{c|cccccccccccccc}
$j$ & $0$ & $1$ & $2$ &   $3$ &   $4$ &   $5$ &  $6$ &  $7$ &  $8$ &  $9$ & $ 10 $ & $ 11 $ & $ \dots $ \\  \hline
$Dh_Z(j)$ & $1$ & $2$ &   $3$ &   $4$ &   $5$ &  $5$ & $5$ & $4$ & $3$ & $2$ & $ 1 $ & $ 0 $ & $ \cdots $ \cr
\end{tabular} \end{center}
and $ \ell(B) = 17 $, $ \ell(Z) = 35 $; \\
\textbf{Case 3}: 
\begin{center}\begin{tabular}{c|cccccccccccccc}
$j$ & $0$ & $1$ & $2$ &   $3$ &   $4$ &   $5$ &  $6$ &  $7$ &  $8$ &  $9$ & $ 10 $ & $ 11 $ & $ \dots $ \\  \hline
$Dh_Z(j)$ & $1$ & $2$ &   $3$ &   $4$ &   $5$ &  $5$ &  $5$ &  $5$ & $3$ & $2$ & $ 1 $ & $ 0 $ & $ \cdots $ \cr
\end{tabular} \end{center}
and $ \ell(B) = 18 $, $ \ell(Z) = 36 $.
\end{prop0}

\begin{proof}
By Proposition \ref{d+1}, we know that $ Dh_{Z}(10) \geq 1 $. Assume that $ Dh_{Z}(10) \geq 2 $. Then, by Proposition \ref{Dav}, necessarily
$ Dh_{Z}(9) \geq Dh_{Z}(10) + 1 \geq 3 $, $ Dh_{Z}(8) \geq Dh_{Z}(9) + 1 \geq 4 $, $ Dh_{Z}(7) \geq Dh_{Z}(8) + 1 \geq 5 $, 
which implies, by \eqref{eq:upbound}, that $  Dh_{Z}(6) \leq 4 $, a contradiction. Therefore
$$ Dh_{Z}(10) = 1. $$
Moreover, for any $ j \geq 11 $ we get  $ Dh_{Z}(j) = 0, $ and also $ Dh_{Z}(9) \geq 2 $, otherwise we violate Proposition \ref{Dav}.\\
Assume that $ Dh_{Z}(9) \geq 3 $. Then, by arguing as before, $ Dh_{Z}(8) \geq Dh_{Z}(9) + 1 \geq 4 $, $ Dh_{Z}(7) \geq Dh_{Z}(8) + 1 \geq 5 $ and $  Dh_{Z}(6) \leq 5 $. If $  Dh_{Z}(6) = 5 $, then $  Dh_{Z}(7) = 5 $ and, by the first line of \eqref{eq:lbound}, $ Dh_{Z}(5) + Dh_{Z}(8) + Dh_{Z}(9) \leq 10 $, so that $ Dh_{Z}(5) \leq 3 $, a contradiction. 
If $ Dh_{Z}(6) \leq 4 $, then, by Proposition \ref{nonincr}, $ Dh_{Z}(7) \leq 4 $, a contradiction. Therefore
$$ Dh_{Z}(9) = 2 $$
and $ Dh_{Z}(8) \geq 3 $, by Proposition \ref{Dav}. \\
Assume that $ Dh_{Z}(8) \geq 4 $. Then, with the same arguments as above, $ Dh_{Z}(7) \geq 5 $. If $ Dh_{Z}(7) = 5 $, then $ Dh_{Z}(6) \geq 5 $ and, being $ \sum_{i = 0}^{10} Dh_{Z}(i) \leq 36 $, then $ Dh_{Z}(5) \leq 4 $, a contradiction; if $ Dh_{Z}(7) > 5 $, then, by\eqref{eq:upbound}, $ Dh_{Z}(6) < 6 $, which provides a contradiction to Proposition \ref{nonincr}. As a consequence, we get that 
$$ Dh_{Z}(8) = 3 $$
so that, by Proposition \ref{Dav} again,  $ Dh_{Z}(7) \geq 4 $. \\
{\bf Case 3.} Assume that $ Dh_{Z}(7) \geq 5 $. If $ Dh_{Z}(7) > 5 $, then, by Proposition \ref{nonincr}, $ Dh_{Z}(6) > 5 $ and so, by \eqref{eq:upbound}, $ Dh_{Z}(6) + Dh_{Z}(7) + Dh_{Z}(8) + Dh_{Z}(9) + Dh_{Z}(10) = 18 $. In particular, it turns to be $ Dh_{Z}(5) = 3 < Dh_{Z}(6)  $, a contradiction. Hence we have
$$ Dh_{Z}(7) = 5. $$
By Proposition \ref{nonincr} we know that $ Dh_{Z}(6) \geq 5 $. In particular, if $ Dh_{Z}(6) > 5 $, then direct computations show that $ Dh_{Z}(5) < 5 $, which violates Proposition \ref{nonincr}. Thus, necessarily,
$$ Dh_{Z}(6) = 5 $$
which implies that 
$$ Dh_{Z}(5) = 5. $$
We get Case 3 of the statement. In particular, we get $ \ell(Z) = 36 $ and $ \ell(B) = 18 $. \\
{\bf Cases 1, 2.} Now assume that 
$$ Dh_{Z}(7) = 4. $$
Therefore, by Proposition \ref{Dav}, we get that $ Dh_{Z}(6) \geq 5 $. If $ Dh_{Z}(6) > 5 $, then, by \eqref{eq:lbound}, $ Dh_{Z}(5) \leq 5 $,  a contradiction. Thus 
$$ Dh_{Z}(6) = 5. $$
As a consequence, by \eqref{eq:lbound} and Proposition \ref{nonincr} we get that either
$$ Dh_{Z}(5) = 6, $$
and we are in Case 1, with $ \ell(Z) = 36 $,  $ \ell(B) = 18 $, or
$$ Dh_{Z}(5) = 5, $$
and we are in Case 2, with  $ \ell(Z) = 35 $, $ \ell(B) = 17 $. 
\end{proof}

We analyze in details the three cases that appear in Proposition \ref{possiblecases}.

\subsection{Case 1} \quad\\
\noindent Let us assume to be in the first case of Proposition \ref{possiblecases}. In fact, we will prove below that this is the general case of the three,
in the sense that Case 2 and Case 3 are limits of Case 1.\\
Recall that the ideal $I_A$ is generated by three quintics $Q_1,Q_2,Q_3$ and one sextic $S$, and its minimal resolution, with Hilbert-Burch matrix $M$,
is described  in Proposition \ref{HBA}.

The Hilbert function tells us that $Z$ is contained in no quintics and it lies in two irreducible sextics. Moreover, we know that $ Z $ satisfies $ CB(9) $,
 and its Hilbert function is the same as the one of a complete intersection of type $ (6,6) $. By the Main Theorem of \cite{Davis84}, 
 it follows that $ Z $ \emph{is} a complete intersection of type $ (6,6) $. So $ I_{Z}$ is generated by $F,F'$, where $ F,F' \in (I_{A})_6 $.
 We get a commutative diagram
 $$
\begin{matrix}
0 & \longrightarrow & R(-12)                   & \xrightarrow  {\begin{pmatrix} F' \\ -F\end{pmatrix}} &  R(-6)^2                              & \longrightarrow  & I_Z                   &\longrightarrow & 0  \\
   &                            & \downarrow M_1 &                                                                                               & \downarrow    M_2             &                            &  \downarrow   &                          &      \\
0 &\longrightarrow  & R(-7)^3                 & \xrightarrow M                                                                      & R(-5)^3\oplus R(-6)            & \longrightarrow  & I_A                  &\longrightarrow & 0
\end{matrix}
$$
where $ M $ is described in  \eqref{M}, and  the matrices $M_1, M_2$ are defined by writing
$$
F = L_1Q_1 + L_2Q_2 + L_3Q_3 + aS \quad F' = L_1'Q_1 + L_2'Q_2 + L_3'Q_3 + a'S
$$
so that 
$$M_{2} = \begin{pmatrix} L _1&  L'_1 \\  L_2 & L'_2 \\ L_3 & L'_3 \\ a & a' \\ \end{pmatrix}, \quad M_{1} = \begin{pmatrix} q_{1} \\ q_{2} \\ q_{3}  \end{pmatrix},$$
with $ q_{i} \in S^{5} \C^{3} $ for $ i\in \{1,2,3\} $. 

By the mapping cone, a resolution of $I_B$ is given by
\begin{equation}\label{idBW}
0 \rightarrow R(-7)^{3} \oplus R(-6)  \xrightarrow {N}  R(-6)^{2} \oplus R(-5)^{3} \rightarrow I_B \rightarrow 0
\end{equation}
where
$$ N = \begin{pmatrix} L_1 & L_{2} & L_{3} & a \\ L'_1& L'_{2} & L'_{3} & a' \\  c_{11} & c_{21} & c_{31}  &  \ell_1 \\ c_{12} & c_{22} & c_{32}  &  \ell_2 \\ c_{13} & c_{23} & c_{33}  &  \ell_3
 \end{pmatrix}. $$ 
and notice that the last three lines represent the transpose of $M$.

\begin{rem0} Let $S'$ be any form of degree $6$ in $I_A$. Then $S'$ corresponds to an element of the homogeneous piece $(I_A)_6$, which is a 
$10$-dimensional linear space. By the description of the generators of $I_A$, there are linear forms $u_1,u_2,u_3$ and a constant $c\in\C$ such that 
$S'$ is the determinant of the matrix:
$$\begin{pmatrix} u_1 & u_{2} & u_{3} & c \\ c_{11} & c_{21} & c_{31} & \ell_{1} \\ c_{12} & c_{22} & c_{32} & \ell_{2}\\ c_{13} & c_{23} & c_{33} & \ell_3.  
\end{pmatrix}$$
If we pose 
$$u_1 = a_{0}x_{0}+a_{1}x_{1}+a_{2}x_{2}, \, u_{2} = a_{3}x_{0}+a_{4}x_{1}+a_{5}x_{2}, \, u_{3} = a_{6}x_{0}+a_{7}x_{1}+a_{8}x_{2}$$
then we can associate to $S'$ the $10$ coordinates $(a_0,\dots,a_8,c)$. This provides a set of coordinates of $S'$ in the linear space  $(I_A)_6$.
\end{rem0}

In particular,  we have
\begin{equation}\label{eq:F1}
F = - \begin{vmatrix} L_1 & L_{2} & L_{3} & a \\ c_{11} & c_{21} & c_{31} & \ell_{1} \\ c_{12} & c_{22} & c_{32} & \ell_{2}\\ c_{13} & c_{23} & c_{33} & \ell_3  \end{vmatrix} 
\end{equation}
\begin{equation}\label{eq:F2}
F' = \begin{vmatrix} L'_1 & L'_{2} & L'_{3} & a' \\ c_{11} & c_{21} & c_{31} & \ell_{1} \\ c_{12} & c_{22} & c_{32} & \ell_{2}\\ c_{13} & c_{23} & c_{33} & \ell_3   \end{vmatrix},
\end{equation} 
so that, if we pose {\small
\begin{multline} \label{eq:lin}
L_1 = a_{0}x_{0}+a_{1}x_{1}+a_{2}x_{2}, \, L_{2} = a_{3}x_{0}+a_{4}x_{1}+a_{5}x_{2}, \, L_{3} = a_{6}x_{0}+a_{7}x_{1}+a_{8}x_{2}  \\
L'_1= a_{10}x_{0}+a_{11}x_{1}+a_{12}x_{2}, \, L'_{2} = a_{13}x_{0}+a_{14}x_{1}+a_{15}x_{2}, \, L'_{3} = a_{16}x_{0}+a_{17}x_{1}+a_{18}x_{2},
\end{multline}
}
then we can associate to the pair $ (F,F') $ the matrix
$$ W = \begin{pmatrix} a_{0} & a_{1} & a_{2} & a_{3} & a_{4} & a_{5} & a_{6} & a_{7} & a_{8} & a_{9} \\
a_{10} & a_{11} & a_{12} & a_{13} & a_{14} & a_{15} & a_{16} & a_{17} & a_{18} & a_{19} \\ \end{pmatrix},$$
where $a_9=a, a_{19}=a'$. The pair $(F, F')$ determines a two dimensional linear subspace $\Lambda$ of $(I_A)_6$, hence an element of the Grassmannian of
$2$-dimensional linear subspaces in a (10)-dimensional space. The Pl\"ucker coordinates of $\Lambda$ are precisely the $2\times 2$ minors of $W$,
in the sense of  \eqref{eq:Pl}.

Since the residue scheme $B$ depends only on the linear space $\Lambda$ spanned by $F,F'$, we can summarize the analysis in the following remark.

\begin{rem0}\label{summa} Let $T\in S^9\C^{3}$ be a form with a decomposition $A\subset \Pj^2$ which satisfies conditions (i)-(iii) above. 
Assume that there exists a second
decomposition $B$ for $T$, such that $Z=A\cup B$ has the difference Hilbert function of Case 1. Then there exists a $2$-dimensional
subspace $\Lambda\subset (I_A)_6$ such that $B$ is linked to $A$ by two sextics which give a basis of $\Lambda$.
\end{rem0}

We can prove that the converse holds. 

\begin{rem0} \label{converse} Choose a general $2$-dimensional subspace $\Lambda\subset (I_A)_6$ and a basis
$F,F'$ of $\Lambda$. Since $I_A$ is generated in degree $6$ and $\Lambda$ is general, the residue $B$ of $A$ in the complete intersection $Z$
of $F,F'$ is a set of $18$ distinct points which does not intersect $A$ (i.e. $Z$ is smooth, as a consequence of the classical Bertini's Theorem). 
Moreover, since $h_Z(9)=35<\ell(Z)=36$, then by the Grassmann formula
the linear spaces $\langle v_9(A)\rangle$ and $\langle v_9(B)\rangle$ meet in exactly one point $T$. Thus, $T$ is a form of degree $9$ 
which has two different (disjoint) decompositions of length $18$.
\end{rem0}

Remark \ref{converse} tells us that we have a (rational) map
\begin{equation}\label{effe}    f : \G(2,10) \dasharrow \langle v_9(A) \rangle \end{equation}
whose image contains tensors $T$ with a second decomposition $B$ of length $18$, such that $Z=A\cup B$ has a Hilbert function as in Case 1 above.
Since the Grassmannian of $2$-dimensional subspaces of $(I_A)_6$ has dimension $2(10-2)=16$, while the span of $v_9(A)$ is $17$-dimensional,
it turns out that a general element  $T\in \langle v_9(A)\rangle$ cannot belong to the (closure of) the image of $f$. 
\smallskip

Next target is to analyze the map $f$ and determine when a given $T$ belongs to the closure of $\imm(f)$.
\medskip

From the Hilbert functions of $A$ and $B$, we know that both $(I_A)_9$ and $(I_B)_9$ determine $37$-dimensional subspaces of the linear space
$R_9$, which has dimension $55$. By standard facts on the intersection of ideals, $I_A\cap I_B$ is the ideal of the union $Z=A\cup B$. Thus, from the Hilbert
function of $Z$ one knows that 
$$ \dim( (I_A)_9\cap (I_B)_9) = \dim((I_Z)_9)= 55-35 = 20.$$
From the Grassmann formula one computes:
$$\dim( (I_A)_9+ (I_B)_9 ) = 37+37-20=54.$$
It follows that $ (I_A)_9+ (I_B)_9$ is a hyperplane in $R_9$,  thus it determines a point in the dual projective space
$\Pj^\vee=\Pj(R_9)^\vee$, of dimension $54$.

Before going on, we need to make a remark on the relation between points and forms in projective spaces, clarifying the roles
of elements of $\Pj^N$ and its dual. The following remark collects standard facts for the
relations between projective geometry and linear algebra.

\begin{rem0}\label{rem:dualnotation} Let $\Pj(V)$ be a projective space and let $\Pj(V)^\vee$ be its dual, defined as the set of hyperplanes in $\Pj(V)$.
 If we fix coordinates in $\Pj(V)$, then there are dual coordinates in $\Pj(V)^\vee$, defined as follows:
for  any  hyperplane $H$ in $\Pj(V)$, the coefficients of an equation $H$ in the fixed coordinates of $\Pj(V)$ are dual coordinates for the point $[H]$ representing $H$ 
in $\Pj(V)^\vee$. Thus, for a point $T\in\Pj(V)$, the coordinates of $T$ are coefficients for an equation of the hyperplane dual to $T$, in the
dual coordinates of $\Pj(V)^\vee$.\\
If $\Lambda$ is a linear subspace of $\Pj(V)$, the dual subspace $\Lambda^\vee\subset\Pj(V)^\vee$ is the set of hyperplanes containing $\Lambda$.
Thus for $\Lambda_1,\Lambda_2\subset\Pj(V)$ subspaces, the intersection $\Lambda_1^\vee\cap\Lambda_2^\vee$ corresponds to the dual of the linear subspace 
$\Lambda_1+\Lambda_2$.

Coming to our situation, up to now we indicated with $\Pj^2$ the projective space over the space of linear forms in the ring $R=\C[x_0,x_1,x_2]$,
i.e. $\Pj^2=\Pj(R_1)$. There is a  natural interpretation of $\Pj(R_9)^\vee=(\Pj^{54})^\vee$ as $\Pj(Sym_9(R_1^\vee))$.
If $T\in R_9$ is a form, the coefficients of $T$ are coordinates of $T$ in the natural frame defined  by monomials. Thus $T$
represents an equation for the hyperplane of $\Pj(R_9)^\vee$ in the dual set of coordinates.\\
In this interpretation, forms in the span $\Lambda$ of powers $L_1^9,\dots, L_r^9$ of linear forms are equations for hyperplanes in 
$\Pj(R_9)^\vee$ associated to points of $\Lambda$. \\
It follows that, if $A=\{L_1,\dots, L_r\}\subset\Pj^2$, then the space $\Pj((I_A)_9)$ has a natural interpretation as the set of hyperplanes
in $\Pj^{54}$ which are dual to points of $\Lambda$.
\end{rem0} 

A consequence of the previous remark is the interpretation of $ (I_A)_9+ (I_B)_9$ expressed in the following:

\begin{prop0}\label{dualsum} There is a natural interpretation of $\Pj(R_9)^\vee$ so that, in the notation above, the hyperplane $ (I_A)_9+ (I_B)_9$
is the dual of the point $T$.
\end{prop0}
\begin{proof} $T$ belongs to the span of the $L_i^9$'s exactly when the hyperplane dual of $T$ contains the intersection of the
hyperplanes dual to the $L_i^9$. In the dual space, the intersection of the
hyperplanes dual to the $L_i^9$ is defined by forms in $ (I_A)_9$. The claim follows.
\end{proof}

%

Next result describes birationally the geometric locus
(hypersurface) of forms of degree $9$ in the span of $\nu_9(A)$ for which a second decomposition exists. 

\begin{thm0}
Let $f : \G(2,10) \dasharrow \langle v_9(A) \rangle$ be the rational map defined in \eqref{effe} and let $ D \subset \G(2,10) $ 
be the divisor of the Grassmannian which defines $ f $. Then $ D = h_{2,10} $. Therefore $ f $ is a linear projection of 
the Pl\"ucker embedding of $ \G(2,10) $.
\end{thm0}
\begin{proof}
Our aim is to show that, if $H$ is a  hyperplane of $\Pj^{54}$ not containing $\langle v_9(A) \rangle$, then 
$ f^{-1}(H)$ is linearly equivalent to the divisor $ h_{2,10}. $
We will do that for a special hyperplane $H$.

In the dual space $(\Pj^{54})^{\vee}$, $H$ corresponds to a point $[H]$. Therefore, we want to describe for which elements $ W \in \G(2,10) $ 
the hyperplane $H_{f(W)}$  associated to the point $ f(W) $  contains the point $ [H] $. \\
As in Remark \ref{rem:dualnotation}, an equation for $H$, whose coefficients are a set of homogeneous coordinates for the dual point $[H]$, 
determines a form of degree $9$ in $R$, which we call $G_H$. 
By construction, 
$$ H_{f(W)} = \Pj((I_{A})_9 + (I_{B(W)})_9) $$
where $ B(W) = \{P_{1}'(W), \dots, P_{18}'(W)\} \subset \Pj^{2} $ is such that $ Z(W) = A \cup B(W) $ 
is the complete intersection of type $ (6,6) $ associated to $ W $. Thus
\begin{multline}\label{eq:iff1}
[H] \in H_{f(W)}\qquad \mbox{ if and only if } \\ \mbox{there exist } G_{A} \in (I_{A})_9,\, G_{B(W)} \in (I_{B(W)})_9 
\mbox{ such that } G_{H} = G_{A} + G_{B(W)}. 
\end{multline}
The condition now becomes purely algebraic.\\
Take a set $Y$ of homogeneous coordinates of the points of $A$, and call $\rho: R\to \C^{18}$ the evaluation map of forms on the set $Y$.
The natural inclusion of $I_{Z(W)}$ in $I_{B(W)}$ gives rise to an exact sequence
\begin{equation}\label{eq:coomA}
0 \to (I_{Z(W)})_9 \longrightarrow (I_{B(W)})_9 \buildrel \rm {\rho} \over  \longrightarrow \C^{18}
\end{equation} 
where the rigthmost map is not surjective.
\begin{claim0}\label{cl:equivalence}
The following are equivalent:
\begin{itemize}
\item[$(i)$] there exist $G_{A} \in (I_A)_9$, $G_{B(W)} \in (I_{B(W)})_9$ such that $G_H = G_{A} + G_{B(W)} $,
\item[$(ii)$] $\rho(G_H)$ belongs to $\rho((I_{B(W)})_9)$. 
\end{itemize}
\end{claim0}
\begin{proof}
If $(i)$ holds, then
$$ \rho(G_H) = \rho(G_{A} + G_{B(W)}) =\rho(G_{A}) + \rho(G_{B(W)})= \rho(G_{B(W)}) $$ 
which implies $(ii)$. \\
On the other side, assume $(ii)$ so that there exists $G_{B(W)} \in (I_{B(W)})_9 $ such that $\rho(G_{B(W)})=\rho(G_H)$.
Therefore $\rho(G_H-G_{B(W)}) = 0$, so that $G_H-G_{B(W)}$ belongs to $(I_A)_9$  and $(i)$ holds. 
\end{proof}
\indent By combining \eqref{eq:iff1} with Claim \ref{cl:equivalence}, our aim turns out to be the following: given a form $ G_H$ of $R_9 $, we want 
to describe the elements $ W \in \G(2,10) $ such that $\rho(G_H)$ belongs to $\rho((I_{B(W)})_9)$. 

\begin{claim0}\label{resZ'}
For any $ W \in \G(2,10) $ set $Z'(W)=Z(W)\setminus \{P_{17},P_{18}\}$. Then for any form $G$ of degree $9$ in $I_{Z'(W)}$, the residue $\rho(G)$
is fixed, up to scalar multiplication.
\end{claim0}
\begin{proof}

We apply Proposition 5.2.10 of \cite{Migliore} (Mapping cone) to the commutative diagram (see Section \ref{sec:link})
$$\begin{CD}
0 @>>> R(-3) @> {\it \begin{pmatrix} -c \\ \ell \end{pmatrix}}  >> R(-1) \oplus R(-2) @>>> I_{\{P_{17},P_{18}\}} @>>> 0\\
@. @AA  A @AA    A @AAA @.  \\
0 @> >>R(-12) @>>{\it \begin{pmatrix} -F' \\ F  \end{pmatrix}} >R(-6)^{\oplus 2} @>>> I_{Z(W)} @>>> 0
\end{CD}$$
where $c,\ell$ are generators for the ideal of $\{P_{17},P_{18}\}$. We get that $ I_{Z'(W)} $ admits a resolution of the form
\begin{equation}\label{eq:idresmc}
0 \to R(-11) \oplus R(-10) \xrightarrow {\begin{pmatrix} \phi_{1} &  \psi_{1}  \\ \phi_{2} &  \psi_{2} \\ -c & \ell \end{pmatrix}} R(-6)^{\oplus 2} \oplus R(-9) \to  I_{Z'(W)}  \to 0
\end{equation}
where $ \phi_{1}, \phi_{2} \in S^{5}\C^{3}, \psi_{1}, \psi_{2} \in S^{4}\C^{3} $ are defined by the central vertical map of the diagram. 
The ideal $ I_{Z'(W)}$ is thus generated by  $F, F', G_{0} $ with
\begin{equation}\label{eq:F1bis}
F = \begin{vmatrix} \phi_{1} & \psi_{1} \\ - c  & \ell \\ \end{vmatrix} = \ell \phi_{1} + c \psi_{1}
\end{equation}
\begin{equation}\label{eq:F2bis}
F' = \begin{vmatrix} \phi_{2} & \psi_{2} \\ - c  & \ell \\ \end{vmatrix} = \ell \phi_{2} + c \psi_{2}
\end{equation}
\begin{equation}\label{eq:G_0}
G_0 = \begin{vmatrix} \phi_{1} & \psi_{1} \\ \phi_{2} & \psi_{2} \\ \end{vmatrix} = \phi_{1} \psi_{2} - \phi_{2} \psi_{1}.  
\end{equation}
Then for any $G\in (I_{Z'(W)})_9$ we get $G=U_1F_1+U_2F_2+qG_0$, $q\in\C$, where $ U_{1}, U_{2} $ are suitable cubics. 
Since $F_1,F_2$ vanish at $Z(W)$, then the residue $\rho(G) $
is a multiple of $\rho(G_0)$.  
\end{proof}

In particular, for $W\in \G(2,10) $, we can find scalars $\alpha_W, \beta_W$ such that the residue $\rho(G)$ is a scalar multiple
of  $(0,\dots,0,\alpha_W,\beta_W)$ for all $G\in (I_{Z'(W)})_9$.
\smallskip

Fix two non-zero scalars $ \alpha, \beta \in \C-\{0\} $. Let $ v_{\alpha,\beta} = (v_{1}, \ldots, v_{18}) \in \C^{18} $ be defined by 
$$ v_{j} = 
\begin{cases}
0  \,\,\,\,\,\,  j \in \{1, \ldots, 16\} \\
\alpha \,\,\,\,\, j=17 \\
\beta \,\,\,\,\, j = 18 \\
\end{cases} $$
Since the evaluation map $\rho$ surjects, we can find $H$ such that the form $G_H$ satisfies $\rho(G_H)= v_{\alpha,\beta} $. We compute $f^{-1}(H)$
for this hyperplane  $H$. Notice that $H$ cannot contain $\langle v_9(A) \rangle$, for $G_H$ does not vanish at $P_{17},P_{18}$.\\
By Claim \ref{cl:equivalence}, $W$ belongs to $f^{-1}(H)$ if and only if there exists $G\in(I_{B(W)})_9$ with $\rho(G)=\rho(G_H)$. Such a form $G$
vanishes at the points $P_1,\dots,P_{16}$, thus it belongs to $(I_{Z'(W)})_9$. We obtain by Claim \ref{resZ'} that $W\in  f^{-1}(H)$ if and only if 
\begin{equation}\label{eq:hyperplane}
\det \begin{pmatrix} \alpha & \beta \\  \alpha_W& \beta_W\\ \end{pmatrix} = 0. 
\end{equation}
Now we compare the two expressions that we have for $F$, from \eqref{eq:F1} and \eqref{eq:F1bis}. We get that 
$$ \ell \phi_{1} + c \psi_{1} =L_1Q_1+L_2Q_2+L_3Q_3+aS,$$
but $Q_1,Q_2,Q_3,S$ vanish at $P_{17},P_{18}$, thus they belong to the ideal spanned by $\ell,c$. In particular, we can write $Q_i=M_ic+N_i\ell$ and
$S=\bar Mc+\bar N\ell$. Since $\ell,c$ give a complete intersection, so they have only trivial syzygies, we conclude that there exists a form $U$ of degree $3$ such that
\begin{gather*}  \phi_1 = Uc + L_1N_1+L_2N_2+L_3N_3+a\bar N \\
		\psi_1 = -U\ell +L_1M_1+L_2M_2+L_3M_3+a\bar M. 
\end{gather*}
Similarly, if we compare the two expressions that we have for $F'$, from \eqref{eq:F2} and \eqref{eq:F2bis}, we get that 
\begin{gather*}  \phi_2  = U'c + L'_1N_1+L'_2N_2+L'_3N_3+a'\bar N \\
		\psi_2   = -U'\ell +L'_1M_1+L'_2M_2+L'_3M_3+a'\bar M. 
\end{gather*}
for some form $U'$ of degree $3$.\\
Now put the previous expressions in the formula \eqref{eq:G_0} for $G_0$. When we compute the residue of $G_0$ at $Z(W)$,
the contribution of $U,U'$ disappears, since $c,\ell$ vanish at the points $P_{17},P_{18}$. Moreover the $M_i,N_i,\bar M,\bar N$ are fixed and do
not depend on the choice of $W$. It follows that $(\alpha_W,\beta_W)$ is an expression in terms of the $2\times 2$ determinants of the coefficients
of the $L_i$'s and $a$ and the coefficients of the $L'_i$'s and $a'$. That is: $(\alpha_W,\beta_W)$ is a linear expression
in the Pl\"ucker coordinates of $W$. This concludes the proof of the theorem.
\end{proof}

Of course, we do not know which projection of the Grassmannian determines the rational map $f$. It is likely that the projection is highly special.

\begin{claim0}\label{thm:bir}
The map $f : \G(2,10) \dasharrow \langle v_9(A) \rangle $ defined in \eqref{effe} is birational onto the image.
\end{claim0}

\begin{proof}
We aim to show that, for some $ P \in \imm(f) \subset \langle v_9(A) \rangle $, the set $ f^{-1}(P) $ is finite and has degree $ 1 $. \\
We proceed via a computational approach in Macaulay2 \cite{Macaulay2} (over a finite field, but then the proof holds
also over $\C$). For a detailed description of our procedure, we refer to the ancillary file \texttt{nonics3.txt}. \\
We start by selecting a finite set $ A = \{P_1,\ldots,P_{18}\} \subset \Pj^2 $, whose elements have random coefficients. 
Then, we construct the Hilbert-Burch matrix of $ A $ and we fix $ 2 $ forms $ F, F' $ of degree $ 6 $ in $ I_{A} $, so that one is not the multiple of the other.
This is equivalent to a choice of $ 6 $ linear forms $ L_{1}, L_{2}, L_{3}, L'_{1}, L'_{2}, L'_{3} $ and $ 2 $ scalars $ a, a' $ 
(since the choice is general, we may assume that $ \partial_{x_{0}} L_{1} = 1, \partial_{x_{0}} L'_{1} = 0, a = 0, a' = 1 $). We get a residual set $ B_{F,F'} $, 
whose ideal admits a resolution as in \eqref{idBW}.
By means of \eqref{effe}, we compute $ f(F,F') $ and we pose $ P = f(F,F') $. Let $ (p_{0}, \ldots, p_{54}) $ be a representative vector for  $ P $. \\
In order to obtain $ f^{-1}(P) $, in the first and second row of the Hilbert-Burch matrix $ N $ of $ B_{F,F'} $ we change $ L_{j} $ (resp. $ L'_{j} $) 
with $ L_{j} = a_{3j-3}x_{0}+a_{3j-2}x_{1}+a_{3j-1}x_{2} $ (resp. with $ L'_{j} = a_{3j+7}x_{0}+ a_{3j+8}x_{1}+a_{3j+9}x_{2} $), for $ j \in \{1,2,3\} $ 
(with $ a_{0} = 1 $ and $ a_{10} = 0 $) and we consider the $ 55\times 54 $ matrix $ M$Fix$'' $ whose columns provide a set of generators  
for $ (I_{A})_{9}+ (I_{B_{F,F'}})_{9} $. 
 Notice that $ M$Fix$''$ is divided in $ 2 $ blocks: the first $ 37 $ columns have integer entries, while in the last $ 17 $ the entries depend linearly 
 on the 2x2 minors of the matrix 
$$ W = \begin{pmatrix}1 & a_{1} & a_{2} & a_{3} & a_{4} & a_{5} & a_{6} & a_{7} & a_{8} & 0 \\
0 & a_{11} & a_{12} & a_{13} & a_{14} & a_{15} & a_{16} & a_{17} & a_{18} & 1 \\ \end{pmatrix}, $$
i.e. on the Pl\"ucker coordinates of $ W $ in the sense of \eqref{eq:Pl}. According to  the previous remark, write $ MFix'' = (A_{1}|A_{2}) $. Therefore 
$$ f^{-1}(P) = \{(a_{1}, \ldots, a_{8}, a_{11},\ldots, a_{18}) \in \C^{16} \, | \, (p_{0}, \ldots, p_{54})\cdot MFix'' = 0_{1 \times 55}\}. $$
Since $ (p_{0}, \ldots, p_{54})\cdot A_{1} = 0_{1 \times 37} $ provide trivial conditions, then 
\begin{equation}\label{eq:sys}
f^{-1}(P) = \{(a_{1}, \ldots, a_{8}, a_{11},\ldots, a_{18}) \in \C^{16} \, | \, (p_{0}, \ldots, p_{54})\cdot A_{2} = 0_{1 \times 17}\}. 
\end{equation}
We notice that the 17 equations  appearing in \eqref{eq:sys} provide a linear system in the  Pl\"ucker coordinates of $ W $. Our  computations in Macaulay2 show that $ f^{-1}(P) $ has dimension $ 0 $ and degree $ 1 $,
which concludes the proof.
\end{proof}

Claim \ref{thm:bir} implies the following:

\begin{rem0}
If $ T \in S^{9}\C^{3} $ of rank $ 18 $ is a general point in $ \imm(f) $, i.e. a general unidentifiable nonic 
of rank $18$, then there exist \emph{exactly two} finite sets computing the rank of $ T $. 

In other words, the space
$ \langle v_9(A) \rangle $ contains a variety $\Theta$, which is the closure of a linear projection of $ \G(2,10) $, whose general points consist of forms in 
$ S^{9}\C^{3} $ of rank $ 18 $, that admit two finite sets computing the rank. 
\end{rem0}


\begin{rem0}\label{rem:crit} 
Given $ T \in S^{9}\C^{3} $ of rank $ 18 $ with a non-redundant finite set $ A = \{P_{1}, \ldots, P_{18}\} \subset \Pj^{2} $ computing it, 
such that $ k_{4}(A) = 15 $ and $ h_{A}(5) = 18 $, by following the proof of \ref{thm:bir}, in principle one can develop a criterion that establishes the uniqueness 
of such an $ A $. 

The corresponding algorithm is explained in more details in Section \ref{sec:alg}.
\end{rem0}

\subsection{Case 2} \quad\\
\noindent Let us assume to be in the second case of Proposition \ref{possiblecases}. This is the unique case in which the new decomposition has less than
$18$ summands, so that the rank of $T$ is $17$, not $18$. In Section \ref{sec:alg} we will provide a software which can exclude that a given form $T$ falls in this case.
\smallskip

Notice that, by Theorem \ref{thm:Davis}, $ Z$  is contained in a plane quintic. Moreover, passing to cohomology in the exact sequence
$$ 0 \rightarrow (I_{Z})_{s}\rightarrow  R(s) \rightarrow  \C^{35} \rightarrow 0  $$
for $ s \in \{5, 7, 12\} $, we get that $ Z $ is contained in a unique quintic $Q$, and there exists a septic $G$ containing $Z$ and not containing $Q$. Since, $ Z $ satisfies $ CB(9) $ and the Hilbert function of $Z$ is the same as the Hilbert function of a complete intersection of type $ (5,7) $, then, by the Main Theorem of \cite{Davis84}, $ Z $ \emph{is} a complete intersection of type $ (5,7) $. In particular, $ I_{Z} = (Q, G) $, with $ Q \in S^{5} \C^{3} $ and $ G \in S^{7} \C^{3} $ and a minimal resolution of $ I_{Z} $ is given by
$$ 0 \rightarrow R(-12) \xrightarrow {\begin{pmatrix} -G \\ Q \\ \end{pmatrix}} R(-5) \oplus R(-7)  \rightarrow I_{Z} \rightarrow 0. $$ 

Again, fix  three quintics $Q_1,Q_2,Q_3$ and one sextic $S$ that generate the ideal of $A$, so that a minimal resolution of $I_A$, with Hilbert-Burch matrix $M$,
is as described  in Proposition \ref{HBA}.

\begin{rem0}\label{ok5ic}
For any choice of the quintic generators $Q_1,Q_2,Q_3$, the ideal $I_A$ coincides with the ideal generated by the $Q_i$'s in degree $7$.\\
Indeed, it follows from the description of the Hilbert-Burch matrix $M$ that, for a general choice of the quintic forms, there are no relations of degree $7$
involving only $Q_1,Q_2,Q_3$. In other words, there are no non-trivial quadrics $q_1,q_2,q_3$ such that $\sum q_iQ_i$ is the zero polynomial.
It follows that the ideal $J$ generated by $Q_1,Q_2,Q_3$ satisfies $\dim J_7=18$, which is exactly the dimension of $(I_A)_7$, as computed from
the Hilbert function.
\end{rem0}

\begin{rem0}\label{17case} In the space $\Pj^{20}=\Pj(R_5)$ which parameterizes quintic forms (up to scalar multiplication) the generators 
$Q_1,Q_2,Q_3$ of $I_A$ determine a plane $\Pj^2=\Pi$. Fix a general quintic $Q$ that belongs to the plane, so that in particular 
$Q$ belongs to $I_A$. Call $J$ the ideal generated by $Q$. The quotient $(I_A)_7/J_7$ is a vector space of dimension  $12$. \\
The variety $E$ of sets $Z$ containing $A$ and complete  intersection of type $(5,7)$ is thus a  $\Pj((I_A)_7/J_7)$-bundle, 
i.e. a $\Pj^{11}$-bundle, over an open subset of $\Pi$. It has dimension $13$.
 \\
Sets $B$ of $17$ points of $\Pj^2$ linked to $A$ by  a complete  intersection of type $(5,7)$ are parameterized by $E$. 
As in Remark \ref{converse}, we get that a general such set $B$ determines one point of the span of $v_9(A)$ with a second decomposition 
of length $17$. Hence  there exists a rational map
 \begin{equation}\label{effe'} f':E \dasharrow \langle v_9(A) \rangle \end{equation}
 such that the closure $\Theta'$  of the image of $f'$ (of dimension at most $13$) is the closure of the locus of forms $T\in \langle v_9(A) \rangle$ for which $A$ is
 non-redundant but the rank is $17$.
\end{rem0}

We will provide in Section \ref{sec:alg} an algorithm which guarantees that a given form $T$ does not lie in the subvariety $\Theta'$, so it has rank $18$.
In order to produce the algorithm, we need a description of the form $T$ which is the intersection of $ \langle v_9(A) \rangle$ and  
$\langle v_9(B) \rangle$, where $B$  is  linked to $A$ by  a complete  intersection $Q\cap G$ of type $(5,7)$.
\smallskip

\begin{rem0}\label{descr17} Fix a quintic $Q\in I_A$, which can be written as 
$$Q=a_1Q_1+a_2Q_2+a_3Q_3,$$
 for a choice of the scalars $a_i$'s. Notice that
if $a_3\neq 0$ (resp. $a_2\neq 0$. $a_1\neq 0$), then the ideal generated by $Q_1,Q_2,Q_3,Q$ coincides with the ideal generated by $Q_1,Q_2,Q$
(resp. $Q_1,Q,Q_3$, $Q,Q_2,Q_3$). By Remark \ref{ok5ic}, a general septic $G\in I_A$ can be written as 
$$G=q_1Q_1+q_2Q_2+q_3Q_3,$$
for a choice of quadrics $q_1,q_2,q_3$. \\
Let $B=B(Q,G)$ be the residue of $A$ in the complete intersection $Q\cap G$. As in case 1, from the Hilbert functions of $A,B$ and $ Z=A\cup B$, we know that
$(I_A)_9$, resp. $(I_B)_9$, determines a $37$-dimensional, resp. a $38$-dimensional, subspace of the linear space
$R_9$, which has dimension $55$. Moreover $\dim( (I_A)_9+ (I_B)_9 ) = 54$, so that
$ (I_A)_9+ (I_B)_9$ is a hyperplane in $R_9$.\\
As in proposition \ref{dualsum}, the hyperplane $ (I_A)_9+ (I_B)_9$ is dual to the point $T$ of intersection between $\langle v_9(A) \rangle$
and $\langle v_9(B) \rangle$.
\end{rem0}

As in case 1, the mapping cone procedure provides an effective way of computing $(I_B)_9$, thus also the sum $ (I_A)_9+ (I_B)_9$.
Namely, the ideal $I_B$ is defined by the minors of the Hilbert-Burch matrix $M$, to which we add the columns that define $Q,G$. We get the matrix:
\begin{equation}\label{M'}
M' = \begin{pmatrix} a_{1} & a_{2} & a_{3} & 0 \\ q_{1} & q_{2} & q_{3} & 0 \\ c_{11} & c_{21} & c_{31} & \ell_{1} \\ c_{12} & c_{22} & c_{32} & \ell_{2} \\ c_{13} & c_{23} & c_{33} & \ell_{3} 
\end{pmatrix}.
\end{equation}
  In particular, $B$ depends on the choice of forms in the first two rows
 $$ \begin{pmatrix}   a_{1} & a_{2} & a_{3} & 0 \\ q_{1} & q_{2} & q_{3} & 0 \end{pmatrix}.$$

 \begin{thm0}\label{bir17}
The map $f' : E \dasharrow \langle v_9(A) \rangle $ defined in \eqref{effe'} is birational onto the image.
\end{thm0}
  \begin{proof}
As in Theorem \ref{thm:bir}, we proceed via a computational approach in Macaulay2 \cite{Macaulay2}. 
For a detailed description of our procedure, we refer to the ancillary file \texttt{nonics1.txt}. \\
Select a finite set $ A = \{P_1,\ldots,P_{18}\} \subset \Pj^2 $, and construct the Hilbert-Burch matrix $M$ of $ A $.
The choice of a quintic form $Q $ and a septic form $G$ in $ I_{A} $, not multiple of $Q$,
is equivalent to a choice of scalars $a_1,a_2,a_3$ and quadrics $q_1,q_2,q_3$.\\
Assume first that $a_1\neq 0$, so that, after rescaling, 
we may assume $a_1=1$.
In this case notice that $Q_1$ is generated by $Q_2,Q_3,Q$. Thus the residue  $B=B(Q,G)$ of $A$ in the complete intersection $Q\cap G$
is also the residual in the complete intersection of $Q\cap G'$ where $G'=(q_2-a_2q_1)Q_2+(q_3-a_3q_1)Q_3=q_2'Q_2+q_3'Q_3$. I.e., we may assume $q_1=0$.
 In particular $B=B(Q,G)$ depends on the choice of $a_2,a_3$ and the $12$ coefficients $ a_4,\dots,a_{15}$ of the quadrics $q_2,q_3$.\\
 Now we can compute how many choices of parameters $a_2,\dots a_{15}$ determine a given hyperplane $P=(I_A)_9+ (I_{B(Q_0,G_0)})_9$. 
 This can be done by Macaulay2. \\
We choose a nonic $P$ corresponding to a hyperplane $ (I_{A})_{9}+ (I_{B(Q_0,G_0)})_{9} $ where $Q_0,G_0$ are determined randomly, 
 with no coefficients $a_1,\dots,a_{15}$ equal to $0$. The coefficients $(p_0,\dots,p_{54})$ of $P$ in  $(\Pj^{54})^\vee$ can be easily computed via
 the mapping cone. Take the $ 55\times 54 $ matrix $N$ whose columns provide a set of generators  
for $ (I_{A})_{9}+ (I_{B(Q,G)})_{9} $, for $Q,G$ general. Then, in the subset of $\Pi$ in which $a_1\neq 0$
\begin{equation}\label{eq:sys'}
{f'}^{-1}(P) = \{(a_2, \ldots, a_{15}) \in \C^{14} \, | \, (p_0, \ldots, p_{54})\cdot N = 0_{1 \times 54}\}. 
\end{equation}
In practice, in order to simplify the computation, observe that in an open set of $E$, we may also assume that one among
$a_4,\dots,a_{15}$ is equal to $1$ . Thus, in the open set,  we get a parametrization of  $E$
with $13$ parameters, corresponding to the dimension of $E$. By varying the coefficients among $a_1,a_2,a_3$ and among  $a_4, \dots, a_{15}$
 which are set equal to $1$, we obtain a complete scan of $f'^{-1}(P)$.\\
Our  computations in Macaulay2 show that, in any case, $ f'^{-1}(P) $ has dimension $ 0 $ and degree 1, i.e. it is a point.
This concludes the proof.
  \end{proof}  

In the rest of the section, we will prove that Case 2 is a degeneration of Case 1. We keep all the previous notation.

\begin{prop0}\label{limiting2}
Let $ P_{0} \in \{G=0\} $ be a general point. Then  $ Z\cup\{P_0\} $ is a limit for $t=0$ of a family $\{ \tilde{Z}_t\}$ (over a small disc), with $ \tilde{Z}_t$  
complete intersection of type $ (6,6) $ containing $A$ for all $t\neq 0$.
\end{prop0}

\begin{proof}
The proof is direct. Let $ \ell_1,\ell_2$ be general generators of $I_{P_{0}}$. Since $P_0\in \{G=0\}$ then there are forms $E_1,E_2$ of 
degree $6$ such that $S=E_1\ell_1+E_2\ell_2$. Since $G$ is irreducible,  for $t\in \C$ general the forms $\ell_1Q+tE_1, \ell_2Q+tE_2$ 
determine a complete intersection $\tilde Z_t $ of type $(6,6)$. By taking the flat limit for $t=0$, the ideal of $\tilde Z_t$ degenerates to
the ideal generated by $\ell_1Q, \ell_2Q, G$, which is the ideal of $Z\cup\{P_0\}$. 

It remains to prove that we can assume $A\subset \tilde Z_t$. Let $A_t$ be the subset of $\tilde Z_t$ which degenerates to $A$ (note that $A_t$ is well defined since
$t$ moves in a small disc). Since the ideal of $A$ is generated in degree $6$, the same holds for the ideal of $A_t$. Moreover $(I_A)_6$
and $(I_{A_t})_6$  have the same dimension. Thus the space of  sextic curves
containing $A_t$ degenerates to the space of sextic curves containing $A$. It follows that the set of complete intersections of type $(6,6)$ containing
$A_t$ degenerates to the set of complete intersections of type $(6,6)$ containing $A$. In particular, $Z$ belongs to the closure
of the set of complete intersections of type $(6,6)$ containing $A$. 
\end{proof}

\begin{prop0}\label{propCase2}
Let $T$ be a form of degree $9$ with decompositions $A,B$ such that $Z=A\cup B$ has Hilbert function as in Case 2. Then there exists a 
family of tensors $T_t$, over a small disc $\Delta$, such that $T_0=T$, $T_t$ belongs to span of   $v_9(A)$ for all $t\in \Delta$ and 
for $t\neq 0$ there exists a second decomposition $B_t$ of $T_t$, which degenerates to $B$, with $A\cup B_t$ complete intersection
of type $(6,6)$.
\end{prop0}
\begin{proof} The claim is an immediate consequence of Proposition \ref{limiting2}, because for $Z_t$ belonging to both Case 1 and Case 2 we have $h^1_9(Z_t)=1$,
thus the spans of $v_9(A)$ and $v_9(B_t)$ intersect in one point, for all $t\in\Delta$.
\end{proof}

It follows that nonic forms with two decompositions $A,B$ such that $A\cup B$ is as in Case 2 are degeneration of
forms with two decompositions $A,B_t$ with $A\cup B_t$ as in Case 1.

\subsection{Case 3} \quad\\
\noindent Let us consider the third case of Proposition \ref{possiblecases}.
We will prove again that the tensor $T$ is the limit of a family tensors $T_t\in v_9(A)$
with two decompositions whose union is complete intersection of two sextic curves.

As in the previous situation, Theorem \ref{thm:Davis} implies that $ Z $ is contained in a plane quintic. Moreover, the evaluation map determines an exact sequence
$$ 0 \rightarrow (I_{Z})_s \rightarrow R(s) \rightarrow \C^{36} \rightarrow 0  $$
for $ s \in \{5, 8, 9,12\} $. We get that $ Z $ is contained in a unique quintic $Q$, and there exists a pencil of plane curves of degree $8$ containing $Z$ 
and not containing $Q$. Fix  two curves $ O_{i} \in S^{8} \C^{3},\ i=1,2 $ such that $Z\subset Q\cap O_1\cap O_2$.

\begin{claim0}\label{cl:opt}
The set $Z$ as above is limit for $t=0$ of a family  $\{Z_t\}$, containing $A$, such that for $t\neq 0$ the set $Z_t$ belongs to an irreducible quintic $Q$ containing $A$.
\end{claim0}
\begin{proof}
Assume that the quintic $Q$ which contains $Z$ is reducible.
In any case, the Hilbert function of $Z$ coincides with the Hilbert polynomial in degree $10$, thus, by the Castelnuovo-Mumford regularity 
Theorem (see \cite{Castelnuovo93}), the ideal of $Z$ is generated in degree $11$. Link $Z$ first with a complete intersection $5,11$, 
then link the residue with a complete intersection $5,7$.
By using twice the mapping cone (see Section \ref{sec:link}, one realizes that the final residue $W$ is complete intersection of $2$ quartics.

There exists a family $\{W_t\}$ of sets of $16$ points, with $W_0=W$, whose generic element is complete intersection of two quartics and it 
is contained in an irreducible quintic curve which contains $A$.

 Indeed, if $X,X'$ are the equations of two quartics whose intersection is $W$,
then $Q=zX+z'X'$, where $z,z'$ are two linear forms which intersect in some point $P\in Q$. Now, move $Q$ in a family of quintics
$\{Q_t\}$, containing $A$, whose general element is irreducible. Move $P$ in a family of points $P_t$ such that $P_t\in Q_t$. The linear forms
$z,z'$ generalize to two families of linear forms $\{z_t\}$, $\{z'_t\}$ such that $z_t,z'_t$ define $P_t$. Then, there are two families of quartics
$\{X_t\}$, $\{X'_t\}$ such that $Q_t=z_tX_t+z'_tX'_t$ and the families specialize to $X,X'$ for $t=0$. So, just take $W_t$ to be the intersection of 
$X_t,X'_t$.

Once the existence of the family $\{W_t\}$ is established, the existence of $\{Z_t\}$ follows immediately by linking back  the general element
of $\{W_t\}$, first with a complete intersection of type $5,7$, and then with a complete intersection of type $5,11$.
\end{proof}

Thus, since we want to find $Z$ as a limit, we may assume that $Q$ is irreducible. The residue $X$ of $Z$ in the complete intersection $Q\cap O_1$ 
is a scheme of length $4$, whose Hilbert function can be computed from the Hilbert functions of $Z$ and the complete intersection (see Section \ref{sec:link}).
It follows that $X$ is contained in a line. Thus $X$ is complete intersection of a line $\ell$ and a quartic $q$. The mapping cone implies thus that 
the homogeneous  ideal of $Z$ is generated by $Q,O_1,O_2$, and its minimal resolution is given by 

\begin{equation}\label{eq:idZvee}
0 \rightarrow R(-12) \oplus R(-9)  \xrightarrow M R(-5) \oplus R(-8)^{\oplus 2} \rightarrow I_Z \rightarrow 0
\end{equation}
where 
\begin{equation}\label{M588}
M = \begin{pmatrix} S & q_{1} \\ q_{2} & \ell_{1}\\ q & \ell\\ \end{pmatrix}
\end{equation}
denotes the Hilbert-Burch matrix of $ I_Z $.

Note that $ S \in S^7\C^{3} $, $ q_{i} \in S^{4}\C^{3} $ for $ i \in \{1,2\} $,  $ \ell_1 \in S^{1}C^{3} $,  and  $ Q, O_{1}, O _{2} $ coincide, 
respectively, with $ (-1)^{i} $ times the minor obtained by leaving out the $ i $-th row of $ M $, $ i \in \{1,2,3\}. $ 
\smallskip

Let $W$ be the residue of $Z$ in the intersection $O_1,O_2$. The fact that $O_1,O_2$ intersect properly, for a general choice of the two forms of degree $8$
 in the ideal of $Z$, follows from the resolution above. By the mapping cone, $W$ is a set of length $28$, complete intersection of $S$ and $q_1$. Thus, the residue
 of $W$ in the intersection  of $q_1$ and $O_1$ is the set $X$.
 
 Move the set $X$ in $O_1$, to a general set $X_t$ of $4$ points, $q_1$ moves to a quartic containing $X_t$. Taking residues, we obtain a family of sets $\{W_t\}$
 of length $28$ in $O_1$ such that $W_0=W$. For $t\neq 0$ the set $W_t$ is the residue of a general set of $4$ points, which is complete intersection of $2$ quadrics.
 Thus the resolution of the ideal of $W_t$ is:
 
$$  0 \rightarrow R(-10)^{\oplus 2}  \rightarrow R(-4) \oplus R(-8)^{\oplus 2} \rightarrow I_{W_t} \rightarrow 0.$$

The ideal of $W_t$ has two minimal generators in degree $8$, one of which is $O_1$ (fixed). 

\begin{claim0}\label{cl:move}
We can move $X$ to $X_t$ so that 
the second generator $O'_t$ of degree $8$ moves to $O_2$ as $t$ goes to $0$.
\end{claim0}
\begin{proof} By the mapping cone, for $t$ general the residue of $X_t$ in a complete intersection $(4,8)$ is contained in a pencil of curves of degree $8$.
The limit of this pencil determines a pencil in the $2$-dimensional space of curves of degree $8$ through $W$. We need to prove that we can choose
the family $\{ X_t\}$ so that the limit contains $O_2$. But this is clear for all the limits contain $O_1$ and the limit changes if we vary the family. 
\end{proof}

Notice that, for $t$ general, the residue of $W_t$ with respect to $O_1\cap O'_t$ is complete intersection of two sextics. It follows that $Z$ is limit of
a family of $36$ points, whose general element is a complete intersection of two sextics.

Collecting all the previous claims, we get:

\begin{prop0} In Case 3, the form $T$ is the limit of a family tensors $T_t\in v_9(A)$
with two decompositions whose union is complete intersection of two sextic curves.
\end{prop0}
\begin{proof} We know that $Z$ is limit of a family $\{Z_t\}$ whose general element is complete intersection of $2$ sextics.
We can conclude that  $Z$ is limit of a family $\{Z_t\}$ as above, \emph{whose general element contains $A$}, because the family
of complete intersections of type $(6,6)$ containing $Z$ is irreducible, and it is the limit of the set of complete intersections
of type $(6,6)$ containing a general set of $18$ points.
\end{proof}

\section{The case of non-empty intersection}\label{sec:intersection}

In this section we assume that $A$ satisfies the genericity conditions (i)-(iii), but we drop the assumption that $A\cap B$ is empty.
We will see that the case can be characterized in terms of the Case 2 of the previous section.

By arguing as in the proof of Claim 4.2 of \cite{AngeC}, we have the following:

\begin{prop0}\label{prop:cap1}
If $A\cap B\neq \emptyset$, then the cardinality of $B$ is $18$ and the intersection contains only one point.
\end{prop0}

\begin{proof}
Assume $ A \cap B = \{P_{1}, \ldots, P_{s}\} $, with $ 1\leq s <18 $. Fix coordinates $T_1,\dots,T_{18}$ for the points of $\nu_9(A)$ and coordinates
$T'_{s+1},\dots,T'_{\ell(B)}$ for the points of $\nu_9(B\setminus A)$. Then there is a choice of scalars such that
\begin{multline*} T = a_{1}T_1+ \ldots + a_{18}T_{18}= \\ =  b_{1}T_1+ \ldots + b_sT_s+ b_{s+1}T'_{s+1}+
\dots b_{\ell(B)}T'_{\ell(B)}. \end{multline*}
Define
$$ T_{0} = (a_{1}-b_{1}) T_1 + \ldots + (a_{s}-b_{s})T_s+ a_{s+1}T_{s+1}+ \ldots + a_{18}T_{18}=  $$
$$ = b_{s+1}T'_{s+1}+ \ldots + b_{\ell(B)}T'_{\ell(B)}. 
\quad\quad\quad\quad\quad\quad\quad\quad\quad\quad\quad\quad\quad\quad\quad\quad\,\,\, $$ 
$ T_{0} $ is an element of $ S^{9} \C^{3} $ admitting two disjoint decompositions: $ A $ and $ B_{0} = B \setminus A $. Notice that, since 
\begin{equation}\label{T0}
 T = T_{0} + b_{1}T'_1+ \ldots + b_{s}T'_s,  
\end{equation}   
necessarily $ B_{0} $ is non-redundant for $ T_{0} $, otherwise $B$ is redundant for $T$. Denote by $ A' \subset A $ a non-redundant decomposition of $ T_{0} $. 
Notice that $ P_{i} \in A' $ for all $ i \in \{s+1,\ldots, 18\} $, otherwise, by \eqref{T0}, $ A $ is a redundant decomposition for $ T $. It turns out that $ T_{0} $ 
has two non-redundant decompositions, $ A' $ and $ B_{0} $, with $ \ell(A') \leq 18 $ and $ \ell(B_{0}) = \ell(B) - s \leq 18 - s \leq \ell(A') $. Since $ A $ 
satisfies properties (ii) and (iii), then, by Remark \ref{maxKrank}, $ k_{4}(A') = \min \{15, \ell(A')\} $ and $ h_{A'}(5) = \ell(A') \leq 18 $, and so, 
by Theorem \ref{range}, $ B_{0} $ cannot exist, unless $\ell(A')=18$, i.e. $A=A'$. Thus also $ A $ is non-redundant for $ T_{0} $.

The tensor $ T_{0} $ has two non-redundant decompositions $A, B_0$, with $\ell(B_0)\leq \ell(B)-s$, and $ A \cap B_{0} = \emptyset $. 
Now assume that either $s\geq 2$ or $\ell(B)<18$, hence $\ell(B_0)\leq 16$. Thus $Z$,
which is also equal to $ A \cup B_{0} $, is a set of at most $34$ points which, from Corollary \ref{CBconseq},  must satisfy the property $\mathit{CB}(9)$.
It follows that:
$$ \ell(Z)\geq \sum_{i=0}^4 Dh_A(i)+Dh_Z(5)+\sum_{i=6}^{10} Dh_Z(i)\geq  2\sum_{i=0}^4 Dh_A(i)+Dh_Z(5) = 30+Dh_Z(5). $$
It follows $Dh_Z(5)\leq 4$. Thus, by Proposition \ref{nonincr}, $Dh_Z(i)\leq 4$ for all $i>4$. Since $Dh_Z(10)>0$, the difference Hilbert function
$Dh_Z$ cannot be strictly decreasing from $5$ to $10$. This contradicts Proposition \ref{Dav}.\\

\end{proof}

\begin{rem0} Thus, if $A\cap B$ is non-empty, there exists a unique point, say $P_1\in A$, which also belongs to $B$. Put $B_0=B\setminus\{P_1\}$.
With the notation of the proof of Proposition \ref{prop:cap1}, we obtain that $T_0$ has two non-redundant decompositions: $A$ and $B_0$, with
$\ell(B_0)=17$. It follows from Proposition \ref{possiblecases} applied to $T_0$ that $Z=A\cup B_0$ has Hilbert function as in Case 2. Thus,  the previous
analysis of Case 2 shows that $A$ and $B_0$ are linked in a complete intersection of type $5,7$, and $T_0$ is the unique point of intersection
of the spans of $\nu_9(A)$ and $\nu_9(B_0)$.\\
Moreover, the tensor $T$ belongs to the line joining $\nu_9(P_1)$ to $T_0$.
\end{rem0}

\begin{rem0} \label{tangency} Let $A,B$ be as in Case 1 of the previous section. We know that $Z=A\cup B$ is the complete intersection of two sextics $F,F'$,
and the spans of $\nu_9(A)$, $\nu_9(B)$ meet in one point. Assume that $F$ moves in a family, so that the limit curves $F_0$ is still irreducible,
and it is tangent to $F$ in one point of $A$, say in $P_1$. Then $B$ moves in a family of finite sets whose limit $B_0$ is linked to $A$ by a complete intersection
$F_0,F'_0$. It is clear from the construction that  $A,B_0$ share the point $P_1$, so that $ \nu_{9}(P_{1}) $ is the unique point of intersection
of the spans of $\nu_9(A)$, $\nu_9(B_0)$. 
Notice that since the ideal of $A$ is generated in degree $6$, for all sextics $F$ containing $A$ one can find infinitely many sextics $F'_0$, containing $A$ and
tangent to $F$ at $P_1$.
\end{rem0}

Now, we are ready to prove that this case can also be seen as a degeneration of the previous Case 1.

\begin{prop0}\label{tplust0}  When $A\cap B\neq \emptyset$, the form $T$ is the limit of a family of tensors $T_t\in v_9(A)$
with two decompositions whose union is complete intersection of two sextic curves.
\end{prop0}
\begin{proof} We know, by Proposition \ref{prop:cap1}, that $A\cap B$ is one point, say $P_1$, moreover $B_0=B\setminus \{P_1\}$
is a set of $17$ points, linked to $A$ by a complete intersection of type $5,7$. Moreover the span of $\nu_9(A)$ and $\nu_9(B_0)$
meet in a form $T_0$ of rank $17$ and $T$ sits in the line joining $\nu_9(P_1)$ and $T_0$.

By Case 2 of the previous section and by Proposition \ref{propCase2}, we know that $T_0$ is a limit of a family $\{T_t\}$, whose
general element $T_t$ is the intersection of the spans of $\nu_9(A)$ and $\nu_9(B_t)$, where $\{B_t\}$ is a family of finite sets
of cardinality $18$, which are linked to $A$ by a complete intersection of two families of sextics $\{F_t\}, \{F'_t\}$. The limits of the two
families $F_0,F'_0$ intersect in a common quintic $Q$, containing $A\cup B_0=Z$. In particular, $F_0=Q\cup L$, where $L$ is a general line. 

 Now, take a family $\{G_t\}$ of sextics through $A$, such that $G_t$ is tangent to $F_t$ for $t$ general (hence for all $t$). Such a family exists because
 the ideal of $A$ is generated in degree $6$. For each $a\in \C$  consider the family $\{H(a)_t\}$, where $H(a)_t=\{G_t +aF'_t\}$. For $t,a$ general the curves
 $H(a)_t$ and $F'_t$ link $A$ to a finite set $B(a)_t$ of cardinality $18$. By the construction of Case 1 of the previous section,
 the generators of the ideal $I_{B(a)_t}$ depend linearly on $a$. Thus  the hyperplane defined by $I_A+I_{B(a)_t}$ determines a pencil of
 hyperplanes, when $a$ moves. For $a=0$ the hyperplane $I_A+I_{B(a)_t}$ corresponds to the intersection $T(a,t)$ of the spans of
 $\nu_9(A)$ and $\nu_9(B(0)_t)$. Since, for $a=0$, $H_t=G_t$ is tangent to $F_t$ at $P_1$, the intersection is $P_1^9$.
 When $a$ goes to infinity and $t$ goes to $0$, the limit of $T(a,t)$ corresponds to the limit of the points defined by the families $\{F_t\}, \{F'_t\}$,
 hence to $T_0$. For any intermediate $a$, we get as a limit for $t=0$ a point corresponding to hyperplanes of the pencil,
 thus points of the line joining $P^9$ and $T_0$.
 
\end{proof}

At the end of the analysis of Case 1,2,3 and the case of non empty intersection, we see that all the forms for which there exists a second
decomposition are limits of forms described in Case 1 of the previous section. Since the image of a Grassmannian in a projection
is Zariski closed, we can summarize the result of the section in the following

\begin{thm0}\label{th9} Let $F$ be a ternary form of degree $9$, with a Waring expression
$$ F= \lambda_1L_1^9+\dots +\lambda_{18} L_{18}^9$$
and call $A$ the set of points induced by $\{L_1,\dots, L_{18}\}$, so that $A$ is a  decomposition of $F$. 
Assume that $A$ satisfies the following genericity properties 
of the beginning of Section \ref{sec:nonics18}:
\begin{itemize}
\item[(i)] $ A $ is non-redundant;
\item[(ii)] $ k_{4}(A)  =15$;
\item[(iii)] $ h_{A}(5) =18$.
\end{itemize}
Then there exists a second decomposition $B$ of length $\leq18$ for $F$ only if $F$ belongs to a fixed hypersurface
$\Theta$ in the span of $\nu_9(A)$, which is the closed image of a birational projection of the Grassmannian
of lines in $\Pj^9$ in its Pl\"ucker embedding.

Therefore, if $F$ does not belong to the hypersurface $\Theta$, then the rank of $F$ is $18$ and $A$ is the
unique decomposition of $F$.
\end{thm0}

\begin{rem0}\label{panforte} ({\it The Panforte Challange}) The bad locus $\Theta$ is thus a hypersurface in a projective space $\Pj^{17}$. 

As in \cite{AngeC}, Remark 4.10, geometrically $\Theta$ is composed of points in which two folds of the secant variety
$Sec_{18}$ to the Veronese variety $\nu_9(\Pj^2)$ cross each other. Thus, the points of $\Theta$ are singular points
 of the secant variety, which is a hypersurface of $\Pj^{54}$. It turns out that the secant variety has a singular locus of codimension $1$.

We do not know an equation, or even the degree of $\Theta$ in $\Pj^{17}$, not even when the points of $A$ are general.
A computer based calculation did not provide an answer, in a reasonable time.

As a challenge, the second author offers a Panforte (traditional cake of Siena) to the first who will
determine an equation for $\Theta$, for a general choice of the decomposition $A$.
\end{rem0}

\section{The algorithm}\label{sec:alg}

This section is devoted to the algorithm we developed based on the criterion  of  minimality explained in Remark \ref{17case}, Remark \ref{descr17}, 
and Theorem \ref{bir17}. 
Our criterion is effective in the sense of \cite{COttVan17b}.

\subsection{The algorithm} \label{algor} Fix a finite set $ A = \{P_{1}, \ldots, P_{18}\} \subset \Pj^{2} $ and a ternary nonic $T$ in the linear span of $ \nu_{9}(A) $, i.e. 
\[
 T = \sum_{i=1}^{18} \lambda_{i}L^{9}_{i} 
\]
where, according to Remark \ref{rem:dualnotation}, $ L_{i} $ is the linear form associated with  $ P_{i} $ in $\Pj^2= \Pj(R_1)$ 
and $ \lambda_{i} \in \C $, for $ i \in \{1, \ldots, 18\}$. 
For any $ i $, denote by $ \vect{v}_{i} $ a representative vector for $ P_{i} $ and by $ (t_{0}, \ldots, t_{54}) $ the coefficients of $ T $ in the standard monomial basis of degree $ 9 $ in $ 3 $ variables. \\ In order to establish that the given $ T $ has exactly rank $ 18 $, proceed as follows. \\ First, check that the tests
 \begin{enumerate}
  \item[1)] \emph{non-redundancy test}: $\dim \langle \nu_9(\vect{v}_1), \ldots, \nu_9(\vect{v}_{18}) \rangle = 18$
  \item[2)] \emph{fourth Kruskal's rank test}: $k_{4}(A) = 15 $
  \item[3)] \emph{fifth Hilbert function test}: $h_{5}(A) = 18 $
 \end{enumerate}
provide positive answers.\\
If so, construct the ideal $ I_{A} $ and its Hilbert-Burch matrix $ M $ as in \eqref{M}. \\ Add to the transpose of $ M $
the rows
$$ \begin{pmatrix} a_{1} & a_{2} & a_{3} & 0 \cr
q_{1} & q_{2} & q_{3} & 0 \cr
\end{pmatrix} $$
so that we get the Hilbert-Burch matrix $ M' $ of a hypothetical non-redundant decomposition $ B $ of $ T $ with $ \ell(B) = 17 $, as in \eqref{M'}. \\ By arguing as in the proof of Theorem \ref{bir17} and by using the same notation, assume that $ a_{1} = 1 $ and $ q_{1} = 0 $ (the cases $ a_{2} = 1 $, $ q_{2} = 0 $ and $ a_{3} = 1 $, $ q_{3} = 0 $ are similar), take the matrix $ N_{1} $ whose columns yield a set of generators for $ (I_{A})_{9} + (I_{B(Q,G)})_{9} $, for $ Q, G $ general and compute 
$$ d_{1} = \max_{4 \leq i \leq 15} \dim \{(a_{2}, \ldots, a_{15}) \in \C^{14} \, | \, (t_{0}, \ldots, t_{54}) \cdot N_{1} = \underline{0} \cap a_{i} = 0 \}. $$
 Equivalently, in the cases $ a_{2} = 1 $, $ q_{2} = 0 $ and $ a_{3} = 1 $, $ q_{3} = 0 $ compute, respectively, 
$$ d_{2} = \max_{4 \leq i \leq 15} \dim \{(a_{1},a_{3} \ldots, a_{15}) \in \C^{14} \, | \, (t_{0}, \ldots, t_{54}) \cdot N_{2} = \underline{0} \cap a_{i} = 0 \} $$
and
$$ d_{3} = \max_{4 \leq i \leq 15} \dim \{(a_{1},a_{2}, \ldots, a_{15}) \in \C^{14} \, | \, (t_{0}, \ldots, t_{54}) \cdot N_{3} = \underline{0} \cap a_{i} = 0 \}. $$
Thus, if the also the next test is successful:
\begin{enumerate}
\item[4)] $ d_{1} = d_{2} = d_{3} = - 1 $, i.e. $ \dim f'^{-1}(T) = -1 $
\end{enumerate}
then $ f'^{-1}(T) $ is empty and so $ T $ has rank $ 18 $.\\


The algorithm has been implemented in Macaulay2, over the finite field $ \Z_{31991} $. The detailed procedure is contained in the ancillary file \texttt{nonics2.txt}. \\

Some examples of 
ternary nonics of rank $ 18 $ and $ 17 $, with a non-redundant decomposition of length $18$, are presented in the following subsection.

\subsection{Examples}\label{exok}
In Macaulay2, we generated a random collection of $ 18 $ points: 
\[ A =
\begin{bmatrix}
 \vect{v}_i
\end{bmatrix}_{i=1}^{18} = {
\begin{bmatrix}
1 & 1 & 1\\
0 & 1 & 2  \\
-1 & 2 & 1 \\
1 & 2 & 3 \\
1 & -2 & 0\\
2 & 1 & 4 \\
4 & 2 & -3 \\
1 & 5 & 1 \\
5 & 2 & 3 \\
6 & 2 & 3 \\
1 & 7 & 7 \\
1 & 7 & 3 \\
6 & 5 & 4 \\
-7 & 2 & 3 \\
3 & 7 & 4 \\
2 & -5 & 1 \\
6 & 3 & -4 \\
-7 & 6 & 6 
\end{bmatrix}.}
\]

By abuse of notation, define $ A=\{P_1,\dots,P_{18}\} \subset \Pj^{2} $, where $ P_i $ is the projective class of $\vect{v}_i$ and denote by $L_i$ the linear form whose coefficients are given by $\vect{v}_i$.

Test 1) and test 3) show, respectively, that $\dim \langle \nu_9(A) \rangle = \operatorname{rank}([\nu_{9}(\vect{v}_i)]_{i=1}^{18}) = 18$ and $ h_{A}(5)= \operatorname{rank}([\nu_{5}(\vect{v}_i)]_{i=1}^{18}) = 18 $, as required.  Moreover, all the $ 816 $ subsets of $ 15 $ columns of $ [\nu_{4}(\vect{v}_i)]_{i=1}^{18} $ are of rank $ 15 $, so that $ k_{4}(A) = 15 $.  

 \smallskip

\paragraph{\textit{A case of rank $18$}}  Let 
$$ T_{2} = \sum_{i=1}^{18} L_{i}^{9} = $$
{\small{$$ = [4283x_{0}^9-14212x_{0}^8x_{1}+2365x_{0}^7x_{1}^2-11335x_{0}^6x_{1}^3+10354x_{0}^5x_{1}^4-7342x_{0}^4x_{1}^5+11432x_{0}^3x_{1}^6+$$
$$-15881x_{0}^2x_{1}^7-10204x_{0}x_{1}^8-663x_{1}^9-10837x_{0}^8x_{2}-6573x_{0}^7x_{1}x_{2}+6070x_{0}^6x_{1}^2x_{2}-12124x_{0}^5x_{1}^3x_{2}+$$
$$+8455x_{0}^4x_{1}^4x_{2}-9097x_{0}^3x_{1}^5x_{2}+200x_{0}^2x_{1}^6x_{2}+11563x_{0}x_{1}^7x_{2}+11173x_{1}^8x_{2}+2810x_{0}^7x_{2}^2+\quad\quad\quad$$
$$+5187x_{0}^6x_{1}x_{2}^2-1688x_{0}^5x_{1}^2x_{2}^2-3089x_{0}^4x_{1}^3x_{2}^2+8745x_{0}^3x_{1}^4x_{2}^2+12508x_{0}^2x_{1}^5x_{2}^2+151x_{0}x_{1}^6x_{2}^2+\quad\quad$$
$$+11119x_{1}^7x_{2}^2+11414x_{0}^6x_{2}^3+2714x_{0}^5x_{1}x_{2}^3+11939x_{0}^4x_{1}^2x_{2}^3+5024x_{0}^3x_{1}^3x_{2}^3+10884x_{0}^2x_{1}^4x_{2}^3+\quad\quad$$
$$+8404x_{0}x_{1}^5x_{2}^3+755x_{1}^6x_{2}^3+15891x_{0}^5x_{2}^4-1013x_{0}^4x_{1}x_{2}^4-11790x_{0}^3x_{1}^2x_{2}^4+14982x_{0}^2x_{1}^3x_{2}^4+\quad\quad\quad$$
$$-8411x_{0}x_{1}^4x_{2}^4-5236x_{1}^5x_{2}^4+4416x_{0}^4x_{2}^5-11481x_{0}^3x_{1}x_{2}^5+14698x_{0}^2x_{1}^2x_{2}^5+5309x_{0}x_{1}^3x_{2}^5+\quad\quad\quad$$
$$+11614x_{1}^4x_{2}^5-9777x_{0}^3x_{2}^6-2702x_{0}^2x_{1}x_{2}^6-5846x_{0}x_{1}^2x_{2}^6-10960x_{1}^3x_{2}^6-8430x_{0}^2x_{2}^7+7085x_{0}x_{1}x_{2}^7+$$
$$+12763x_{1}^2x_{2}^7-14136x_{0}x_{2}^8-9808x_{1}x_{2}^8+9194x_{2}^9].\quad\quad\quad\quad\quad\quad\quad\quad\quad\quad\quad\quad\quad\quad\quad\quad\quad\quad$$}}

\noindent Since tests 1), 2) and 3) are successful, then $ T_{2} $ is general enough, so that our criterion applies. Moreover, test 4)  provides positive answer too. Therefore we conclude that $ A $ is minimal for $ T_{2} $, i.e. $T_2$ has rank $18$ (ancillary \texttt{nonics2.txt}). \\

\paragraph{\textit{A case of lower rank}} In the same span of $v_9(A)$ one can find forms for which the decomposition $A$ is non-redundant, yet there is another decomposition of length $17$.

For instance, one can take {\small{$$ (\lambda_{1},\ldots, \lambda_{18}) = \\
(10308,-9437,-13956,-12270,2135, -4854,-2213,1755,-13629,$$
$$\quad\quad\quad\quad\quad 7308,-8496, 2940,11348,-12437,-6712,4086,-823,-2818)$$}}
so that
$$ T_{1} = \sum_{i=1}^{18}\lambda_{i} L_{i}^{9} = $$
{\small{$$ =[9666x_0^9+13004x_0^8x_1+12463x_0^7x_1^2-13235x_0^6x_1^3 -15442x_0^5x_1^4+15509x_0^4x_1^5+
 -6311x_0^3x_1^6+$$
$$ \,\,\,\,\,\,\, -2390x_0^2x_1^7+547x_0x_1^8-119x_1^9-14916x_0^8x_2+1822x_0^7x_1x_2-8022x_0^6x_1^2x_2-9386x_0^5x_1^3x_2+$$
$$\,\,\,\,\,\,\,-2742x_0^4x_1^4x_2 +10541x_0^3x_1^5x_2 +1156x_0^2x_1^6x_2 -12023x_0x_1^7x_2
+ 4417 x_1^8x_2-11823x_0^7x_3^2 -737x_0^6x_1x_1^2+$$
$$\,\,\,\,\,\,\,-7616x_0^5x_1^2x_2^2+11293 x_0^4x_1^3x_2^2 -8260x_0^3x_1^4x_2^2-9332x_0^2x_1^5x_2^2 +7078x_0x_1^6x_2^2-4553x_1^7x_2^2-15941x_0^6x_2^3+$$
$$  +4339x_0^5x_1x_2^3-4251x_0^4x_1^2x_2^3+ 9854x_0^3x_1^3x_2^3 -22x_0^2x_1^4x_2^3 + 8408x_0x_1^5x_2^3 +11858x_1^6x_2^3+ $$
 $$ \,\,\,\,\,\,\,-9161x_0^5x_2^4-9854x_0^4x_1x_2^4 -13165x_0^3x_1^2x_2^4 -2105x_0^2x_1^3x_2^4 -8715x_0x_1^4x_2^4 +390x_1^5x_2^4 -9955x_0^4x_2^5+$$
$$\,\,\,\,\,\,\, -11013x_0^3x_1x_2^5 -10651x_0^2x_1^2x_2^5 -3850x_0x_1^3x_2^5+4029x_1^4x_2^5-11735x_0^3x_2^6 -12427x_0^2x_1x_2^6+ 12255x_0x_1^2x_2^6 +$$
$$\,\,\,\,-3686x_1^3x_2^6-2271x_0^2x_2^7+5939x_0x_1x_2^7-3402x_1^2x_2^7 +13298x_0x_2^8+6455x_1x_2^8+x_2^9].\quad\quad  $$}}

\noindent $T_1$ has been obtained as the point of intersection of the span of $\nu_9(A)$ with $\nu_9(B)$, where $B$ is a set of $17$ points, linked to $A$ with two curves $ Q_{0} $ and $ G_{0} $ of degree $ 5 $ and $ 7 $ in $ I_{A} $. In particular, $ Q $ (resp. $ G_{0} $) is defined by the determinant of the $4 \times 4$ matrix obtained by adding the row $ R_{1} = (1,10399,13534,0) $ (resp. the row $ R_{2} = (0,-633x_0^2-11455x_0x_1+2134x_0x_2+11038x_1^2-8888x_1x_2-588x_2^2, 1927x_0^2+4114x_0x_1+11328x_0x_2+13814x_1^2-10664x_1x_2-1749x_2^2,0) $) to the transpose of the Hilbert-Burch matrix $ M $ of $ A $, i.e.
$$ (a_{1},a_{2}, \ldots, a_{15}) = \\
(1,10399,13534,-633,-11455,2134,11038,-8888,$$
$$ \quad\quad\quad\quad\quad\quad\,\,\,-588,1927,4114,11328,13814,-10664,-1749).$$
 Notice that, according to the notation of test 4), in this case $ d_{1} = d_{2} = d_{3} = 0 $ and 
$$ \max_{4 \leq i \leq 15} \deg \{(a_{2}, \ldots, a_{15}) \in \C^{14} \, | \, (t_{0}, \ldots, t_{54}) \cdot N_{1} = \underline{0} \cap a_{i} = 0 \} = 0, $$
$$ \max_{4 \leq i \leq 15} \deg \{(a_{1},a_{3} \ldots, a_{15}) \in \C^{14} \, | \, (t_{0}, \ldots, t_{54}) \cdot N_{2} = \underline{0} \cap a_{i} = 0 \} = 0, $$
$$ \max_{4 \leq i \leq 15} \deg \{(a_{1},a_{2}, \ldots, a_{15}) \in \C^{14} \, | \, (t_{0}, \ldots, t_{54}) \cdot N_{3} = \underline{0} \cap a_{i} = 0 \} = 0. $$
Therefore, the rank of $ T_1 $ is at most $ 17 $. Indeed, the rank is exactly $17$, by Proposition \ref{bir17} and its proof (ancillary file \texttt{nonics1.txt}).\\
In order to get coordinates for the points in $ B $, one needs to solve the polynomial system given by the maximal minors of the $ 5 \times 4 $ matrix $ {\small{\begin{pmatrix} R_{1} \cr R_{2} \cr M^{t} \\ \end{pmatrix}}} $. This can be achieved with Macaulay2 software system: indeed, by computing the eigenvalues and eigenvectors of certain \emph{companion} matrices (ancillary file \texttt{nonics4.txt}) one can find the following representative vectors in $ \C^2 $ for the points in $ B $
{\small{$$ (1,62.6659,29.7378) $$
$$ (1,13.368+38.1825 \, i,-19.099+7.53788 \, i) $$
$$ (1,13.368-38.1825 \, i,-19.099-7.53788 \, i) $$
$$ (1,35.333,40.797) $$
$$ (1,14.7061,27.8538) $$
$$ (1,10.7119,4.95399) $$
$$ (1,-0.796312,2.23381) $$
$$ (1,1.06064+0.13583 \, i,1.62951-0.563286 \, i) $$
$$ (1,1.06064 -0.13583 \,i,1.62951+0.563286 \, i) $$
$$ (1,0.737271,-0.0631582) $$  
$$ (1,-0.245331,-0.76262) $$
$$ (1,-0.187307,0.100519) $$
$$ (1,-0.0870499,-0.126324) $$
$$ (1,0.00104432,0.00164595) $$
$$ (1,0.306581+0.0182712 \, i,-0.877193-0.031211\, i) $$
$$ (1,0.306581-0.0182712 \, i,-0.877193+0.031211 \, i) $$
$$ (1,0.390447,0.585521). $$}}

\begin{rem0}
In principle one can try to use our analysis also to determine the identifiability of a form $T$. In this case one starts by adding, in the previous algorithm,
the following test:
 \begin{enumerate}
\item[5)] check that the solution set of the polynomial system $(t_{0}, \ldots, t_{54})\cdot A_{2} = 0_{1 \times 17}$ introduced in \eqref{eq:sys}
has dimension $-1$. 
 \end{enumerate}
If the answer is negative, we can conclude that there are other solutions of the system \eqref{eq:sys}, so one can guess that $T$ has 
other decompositions.\\
Unfortunately, if the answer is positive, one cannot immediately conclude the identifiability of $T$. Namely, test 5) checks that there exists no other
set $B$ of length $18$, linked to $A$ by a complete intersection of two sextics $F,F'$, such that $T$ also sits in the span of $\nu_9(B)$,
 {\it but only when $A\cap B=\emptyset$ and $F,F'$ have no common components}.\\
However, we know that there are limit cases (the case $A\cap B$ non-empty, or the case 3 of Proposition \ref{possiblecases}) in which
$B$ exists but the intersection $F\cap F'$ is not proper. To exclude these cases one needs new ad hoc tests, which can be constructed in principle,
but then the procedure becomes quite laborious.\\
We believe that the natural  way to prove identifiability is to produce equations for the locus $\Theta$ described in remark \ref{panforte}, and test the vanishing 
of the equation for $T$.
\end{rem0}

\paragraph{\textit{An unidentifiable case}}  
Of course, one can use the construction of Case 1 of Proposition \ref{possiblecases} to produce many examples 
of unidentifiable forms.\\
Let {\small{$$ (\lambda_{1}, \ldots, \lambda_{18}) = (5864,9496,11539,1233,-13315,-14222,10709, -5067, 13797,$$
$$ \quad\quad\quad\quad\quad\quad\,\,\,\, 13169,-10531,1592,12589,1728,-4725,-4784,-8696,7515) $$}}
and let
$$ T_{3} = \sum_{i=1}^{18}a_{i} L_{i}^{9} = $$
{\small{$$ = [11096x_{0}^9+14876x_{0}^8x_{1}+14398x_{0}^7x_{1}^2-11088x_{0}^6x_{1}^3-3441x_{0}^5x_{1}^4+11138x_{0}^4x_{1}^5-3819x_{0}^3x_{1}^6+$$
$$+6626x_{0}^2x_{1}^7-9525x_{0}x_{1}^8-9028x_{1}^9+11951x_{0}^8x_{2}-14433x_{0}^7x_{1}x_{2}+15878x_{0}^6x_{1}^2x_{2}+3683x_{0}^5x_{1}^3x_{2}+$$
$$+12902x_{0}^4x_{1}^4x_{2}+9968x_{0}^3x_{1}^5x_{2}+1167x_{0}^2x_{1}^6x_{2}-1011x_{0}x_{1}^7x_{2}+11114x_{1}^8x_{2}-1174x_{0}^7x_{2}^2+\quad\quad$$
$$-10039x_{0}^6x_{1}x_{2}^2+15571x_{0}^5x_{1}^2x_{2}^2-1797x_{0}^4x_{1}^3x_{2}^2+7799x_{0}^3x_{1}^4x_{2}^2+3353x_{0}^2x_{1}^5x_{2}^2-9008x_{0}x_{1}^6x_{2}^2+$$
$$+7892x_{1}^7x_{2}^2+8863x_{0}^6x_{2}^3+12538x_{0}^5x_{1}x_{2}3+4156x_{0}^4x_{1}^2x_{2}^3+5014x_{0}^3x_{1}^3x_{2}^3+15217x_{0}^2x_{1}^4x_{2}^3+\quad\,\,\,$$
$$+10693x_{0}x_{1}^5x_{2}^3-14254x_{1}^6x_{2}^3-12480x_{0}^5x_{2}^4+15094x_{0}^4x_{1}x_{2}^4+11796x_{0}^3x_{1}^2x_{2}^4-11496x_{0}^2x_{1}^3x_{2}^4+$$
$$-3087x_{0}x_{1}^4x_{2}^4-7767x_{1}^5x_{2}^4+1751x_{0}^4x_{2}^5+9059x_{0}^3x_{1}x_{2}^5+14238x_{0}^2x_{1}^2x_{2}^5-640x_{0}x_{1}^3x_{2}^5+\quad\quad\quad$$
$$-14792x_{1}^4x_{2}^5-14262x_{0}^3x_{2}^6-6895x_{0}^2x_{1}x_{2}^6+13550x_{0}x_{1}^2x_{2}^6+7631x_{1}^3x_{2}^6+9523x_{0}^2x_{2}^7+\quad\quad\quad\quad$$
$$-2161x_{0}x_{1}x_{2}^7-3449x_{1}^2x_{2}^7-7220x_{0}x_{2}^8+395x_{1}x_{2}^8+x_{2}^9].\quad\quad\quad\quad\quad\quad \quad\quad\quad\quad\quad\quad \quad\quad\quad   $$}}

\noindent As in the previous cases, tests 1), 2), and 3) are successful for $ T_{3} $. On the other hand, in this case our computations show that the polynomial system $(t_{0}, \ldots, t_{54})\cdot A_{2} = 0_{1 \times 17}$ defined in \eqref{eq:sys} has a solution set of dimension $0$ and degree $1$. 
In particular, $ T_{3} $ admits at least two decompositions of length $18$, $ A $ and another set $B$.\\
\noindent Indeed, $ T_{3} = f(F,F') $, where $ F $ (resp. $ F' $) is the ternary form of degree $ 6 $ defined by the determinant of the $ 4\times 4 $ matrix obtained by adding the row $ (x_{0}+14307x_{1}+13416x_{2}, 11657x_{0}+ 9248x_{1}+ 8324x_{2}, -13193x_{0}-1403x_{1}+12171x_{2}, 0) $ (resp. the row $(7694x_{1}+12549x_{2}, -12983x_{0}+538x_{1}+11728x_{2}, 743x_{0}-12966x_{1}+12870x_{2},1)$) to the transpose of the Hilbert-Burch matrix $ M $ of $ A $. Therefore, $ T_{3} $ is computed by two non-redundant finite sets of length $ 18 $: $ A $ and its residual set, $B= B_{F,F'} $, in the complete intersection $ (6,6) $ given by $ F $ and $ F' $. According to our theory, test 5) fails for $ T_{2} $. Notice that test 4) is successful for $ T_{3} $, which means that $ T_{3} $ has rank $ 18 $.  CITARE ancillary nonics3.txt con referenza arXiv

\noindent A non trivial solution for the previous system determines the Hilbert-Burch matrix for $B$. Then, coordinates for the points of $B$ can be found by solving (e.g. with the software Macaulay2) the polynomial system given by the maximal minors of its Hilbert-Burch matrix. Indeed, by computing the eigenvalues and eigenvectors of certain \emph{companion} matrices (for more details on the procedure see the ancillary file {\tt{nonics5.txt}}), we get that the following representative vectors in $ \C^{3} $ for the points of $ B $:
 
{\small{$$ (1,-7.96881+4.74847 \,i ,29.737-8.31447+5.39065 \,i) $$
$$ (1,-7.96881-4.74847\, i,-8.31447-5.39065\, i) $$
$$ (1,-8.88473,-21.3598) $$
$$ (1,-3.19251,-1.98613) $$
$$ (1,2.29572+0.339361 \, i,2.08576+2.10835\, i) $$
$$ (1,2.29572-0.339361 \, i,2.08576-2.10835 \, i) $$
$$ (1,1.1725,0.789914) $$
$$ (1,-0.662147+0.268568 \, i,1.41128+0.060661 \, i) $$
$$ (1,-0.662147-0.268568 \, i,1.41128-0.060661 \, i) $$
$$ (1,-0.676455+0.162048 \,i,0.269336-0.242414 \, i) $$
$$ (1,-0.676455-0.162048 \, i, 0.269336+0.242414 \, i) $$
$$ (1,0.299266+0.586034 \, i,0.441543+0.153418 \, i) $$
$$ (1,0.299266-0.586034 \, i, 0.441543-0.153418 \, i) $$
$$ (1,-0.0000365176,-0.00000906466) $$
$$ (1,0.209511,0.479921) $$
$$ (1,0.483517+0.0585949 \, i,-0.520821-0.078578\, i) $$
$$ (1,0.483517-0.0585949 \, i,-0.520821+0.078578 \, i) $$
$$ (1,0.503511, 0.553533) $$}}.

\bibliographystyle{amsplain}
\bibliography{biblioLuca}

\end{document}